\documentclass{amsart}

\newtheorem{theorem}{Theorem}[section]
\newtheorem{lemma}[theorem]{Lemma}

\theoremstyle{definition}

\theoremstyle{remark}
\newtheorem{remark}[theorem]{Remark}
\theoremstyle{corollary}
\newtheorem{corollary}[theorem]{Corollary}

\usepackage{amsmath}
\usepackage{amssymb}
\usepackage{amstext}
\usepackage{amscd}
\usepackage{arydshln}
\usepackage{color}
\usepackage{graphicx}
\usepackage{longtable}
\usepackage{calc}
\usepackage{color}
\usepackage{microtype}
\usepackage{verbatim}
\graphicspath{{Figures/}}

\newcommand{\diag}[1]{\mbox{${\rm diag}(#1)$}}

\newcommand{\wh}[1]{\widehat{#1}}

\def\e{{\epsilon}}

\newcommand{\beq}{\begin{equation}}
\newcommand{\eeq}{\end{equation}}
\newcommand{\bac}[1]{\left[\begin{array}{#1}}
	\newcommand{\eac}{\end{array}\right]}
\newcommand{\bmx}{\begin{bmatrix}}
	\newcommand{\emx}{\end{bmatrix}}

\newcommand{\la}{\lambda}

\newcommand{\hide}[1]{}

\numberwithin{equation}{section}

\newcounter{parenumi}
	{\end{list}}

\newlength{\tablecolwidth}
\begin{document}

\title[Structural backward stability in rational eigenvalue problems]{Structural backward stability in rational eigenvalue problems solved via block Kronecker linearizations}


\author{Froil\'{a}n M. Dopico}
\address{Departamento de Matem\'aticas, Universidad Carlos III de Madrid, Avda. Universidad 30, 28911 Legan\'es, Madrid, Spain}
\email{dopico@math.uc3m.es}
\thanks{The first and second authors were partially supported by ``Ministerio de Econom\'ia, Industria y Competitividad (MINECO)" of Spain and ``Fondo
	Europeo de Desarrollo Regional (FEDER)" of EU through grant MTM2015-65798-P,
	by the ``Proyecto financiado por la Agencia Estatal de Investigaci\'on de Espa\~na'' (PID2019-106362GB-I00 / AEI / 10.13039/501100011033) and by the Madrid Government (Comunidad de Madrid-Spain) under the ``Multiannual Agreement with Universidad Carlos III de Madrid in the line of Excellence of University Professors (EPUC3M23), and in the context of the V PRICIT (Regional Programme of Research and Technological Innovation)". }

\author{Mar\'{i}a C. Quintana}
\address{Department of Mathematics and Systems Analysis, Aalto University, Otakaari 1, Espoo, Finland}
\email{maria.quintanaponce@aalto.fi}
\thanks{The second author was funded by the “contrato predoctoral” BES-2016-076744 of MINECO and by an Academy of Finland grant (Suomen Akatemian päätös 331240).}

\author{Paul~Van~Dooren}
\address{Department of Mathematical Engineering, Universit\'{e}
	catholique de Louvain, Avenue Georges Lema\^{i}tre 4, B-1348
	Louvain-la-Neuve, Belgium}
\email{paul.vandooren@uclouvain.be}
\thanks{This work was developed while the third author held a ``Chair of Excellence UC3M - Banco de Santander'' at Universidad Carlos III de Madrid in the academic year 2019-2020}

\subjclass[2020]{65F15, 15A18, 15A22, 15A54, 93B18, 93B20, 93B60}

\keywords{rational matrix, rational eigenvalue problem, linearization,  matrix pencils, perturbations, backward error analysis}

\date{}

\dedicatory{}

\begin{abstract}
	In this paper we study the backward stability of running a backward stable eigenstructure solver on a pencil $S(\la)$ that is a strong linearization of a rational matrix $R(\la)$ expressed in the form $R(\la)=D(\la)+ C(\la I_\ell-A)^{-1}B$, where $D(\la)$ is a polynomial matrix and $C(\la I_\ell-A)^{-1}B$ is a minimal state-space realization. We consider the family of block Kronecker linearizations of $R(\la)$, which have the following structure
	$$  S(\la):=\left[\begin{array}{ccc} M(\la) & \widehat K_2^T C &  K_2^T(\la)  \\ B \widehat K_1 & A- \la I_\ell & 0\\
	K_1(\la) &   0 &  0
	\end{array}\right],$$
	where the blocks have some specific structures. Backward stable eigenstructure solvers, such as the  $QZ$ or the staircase algorithms, applied to $S(\la)$ will compute the exact eigenstructure of a perturbed pencil $\wh S(\la):=S(\la)+\Delta_S(\la)$ and the special structure of $S(\la)$ will be lost, including the zero blocks below the anti-diagonal.
	In order to link this perturbed pencil with a nearby rational matrix, we construct in this paper a strictly equivalent pencil $\widetilde S(\la)=(I-X)\wh S(\la)(I-Y)$ that restores the original structure, and hence is a block Kronecker linearization of a perturbed rational matrix  $\widetilde R(\la) = \widetilde D(\la)+ \widetilde C(\la I_\ell- \widetilde A)^{-1} \widetilde B$, where $\widetilde D(\la)$ is a polynomial matrix with the same degree as $D(\la)$. Moreover, we bound appropriate norms of $\widetilde D(\la)- D(\la)$, $\widetilde C - C$, $\widetilde A - A$ and $\widetilde B - B$ in terms of an appropriate norm of $\Delta_S(\la)$. These bounds may be, in general, inadmissibly large, but we also introduce a scaling that allows us to make them satisfactorily tiny, by making the matrices appearing in both $S(\la)$ and $R(\la)$ have norms bounded by $1$. Thus, for this scaled representation, we prove that the staircase and the $QZ$ algorithms compute the exact eigenstructure of a rational matrix $\widetilde R(\la)$ that can be expressed in exactly the same form as $R(\la)$ with the parameters defining the representation very near to those of $R(\la)$. This shows that this approach is backward stable in a structured sense.
\end{abstract}

\maketitle

%
%


\section{Introduction}\label{introduction}
It has been known since the 1970's that the zeros of a rational matrix are also the eigenvalues of an appropriately defined pencil of matrices, i.e., a  polynomial matrix of degree at most $1$, and that its poles are the eigenvalues of a principal submatrix of such a pencil. This connection was established in the influential book of Rosenbrock \cite{Ros70}. About 10 years later numerical algorithms were proposed in \cite{VanDooren79, VanDooren81} to construct such a pencil in a numerically stable way. Not only the zeros and poles can be determined via these pencils, but also their structural indices, or partial multiplicities, as well as the minimal indices of the left and right null-spaces of the rational matrix, see e.g. \cite{VVK79}. Together, these are called the eigenstructure of the rational matrix, and the pencils considered in \cite{VVK79} are called {\em system matrices} of a {\em strongly irreducible generalized state-space realization}.

Polynomial matrices can be viewed as special cases of rational matrices, which happen to have all their poles at infinity. The notions of generalized state-space realizations and corresponding (strongly) irreducible system matrices therefore apply to polynomial matrices as well. But in the classic reference \cite{GLR82} a new notion of {\em strong linearization} is introduced for polynomial matrices which is consistent with that of {\em strongly irreducible system matrix} of \cite{VVK79} for the finite eigenvalues and their structural indices. But for the structural indices at infinity, these two definitions differ by a constant shift, which means that the
structural indices at infinity can easily be recovered from one definition to the other. Moreover, the definition of strong linearization introduced in \cite{GLR82} does not guarantee any relationship between the minimal indices of the linearization and those of the polynomial matrix \cite{dtedopmack2014}, in contrast with the pencils in \cite{VVK79} for which the minimal indices are equal.

Even though the definition of strong linearization in \cite{GLR82} was originally given for polynomial matrices, there have been several attempts to extend it to rational matrices \cite{ADMZ, DasAlam19}, including these extensions also the concept of (non-strong) linearization \cite{AlBe16,ADMZ}. Thus, inspired by previous results for polynomial matrices \cite{DLPV}, a wide family of strong linearizations called {\em strong block minimal bases linearizations} is proposed in \cite[Theorem 5.11]{ADMZ} for any $m\times n$  rational matrix $R(\la)$ with coefficient matrices in an arbitrary field $\mathbb{F}$. These linearizations are based on the splitting of $R(\la)$ into its strictly proper part $R_p(\la)$ and its polynomial part $D(\la)$ and in the representation~:
\begin{equation}\label{eq_rationalmatrix}
R(\la) := R_p(\la) +D(\la) = C(\la I_\ell-A)^{-1}B + \sum_{i=0}^d D_i\la^i,
\end{equation}
where $C(\la I_\ell-A)^{-1}B$  is a minimal state-space realization of the strictly proper part $R_p(\la)$,  represented in what follows by the triple $\{A,B,C\}$, and $d >1$ is the degree of the polynomial part. Then $R(\la)$ is represented by the quadruple $\{\la I_\ell-A,B,C,D(\la)\}$. Since in this paper we are analyzing perturbations related to backward errors of eigenvalue solvers
of pencils with real or complex matrix coefficients, we restrict $\mathbb{F}$ to be the real field $\mathbb{R}$ or the complex field $\mathbb{C}$.

A particular case of the strong block minimal bases linearizations in \cite[Theorem 5.11]{ADMZ}
of any $m\times n$ rational matrix $R(\la)$ represented as in \eqref{eq_rationalmatrix} are (modulo block permutations) the pencils of the form
\begin{equation}\label{block_kronecker_lin}
S(\la):=\left[\begin{array}{ccc} M(\la) & \widehat K_2^T C &  K_2^T(\la)  \\ B \widehat K_1 & A- \la I_\ell  & 0\\
K_1(\la) &   0 &  0
\end{array}\right],
\end{equation}
with
$$
K_1(\la):= L_\epsilon(\la)\otimes I_n,\quad \widehat K_1:=\mathbf{e}_{\epsilon+1}^T\otimes I_n, \quad K_2(\la):= L_\eta(\la)\otimes I_m, \quad \widehat K_2:=\mathbf{e}_{\eta+1}^T\otimes I_m,
$$
and where $\otimes$ denotes the Kronecker product, $\mathbf{e}_k = [0 \cdots 0 \, 1]^T$ is the standard $k$th unit vector of dimension $k$ and $L_k(\la)$ is the classical Kronecker block of dimension $k\times(k+1)$
$$ L_k(\la):=\left[\begin{array}{ccccc} 1 & -\la  \\ & 1 & -\la \\ & & \ddots & \ddots \\ & & & 1 & -\la
\end{array}\right].
$$
Moreover, the block $M(\la)$ in \eqref{block_kronecker_lin} is related to the polynomial part $D(\la)$ in \eqref{eq_rationalmatrix} by the ``dual basis'' vector $ \Lambda_k(\la)$ of powers of $\la$,
$$ \Lambda_k^T(\la) :=\left[\begin{array}{ccccc} \la^k & \cdots  & \la^2 & \la & 1 \end{array}\right],
$$
which satisfies $L_k(\la)\Lambda_k(\la)=0$  and also
$$  D(\la)=(\Lambda_\eta(\la)\otimes I_m)^T M(\la) (\Lambda_\epsilon(\la)\otimes I_n).
$$
Thus, $d = \epsilon + \eta +1$ (see \cite[eq. (4.5)]{DLPV}).
The strong linearizations \eqref{block_kronecker_lin} are inspired by the so-called ``block Kronecker linearizations'' that were introduced in \cite[Section 4]{DLPV} for an arbitrary $m\times n$  polynomial matrix $D(\la)$. Therefore, we use the same name in the rational setting. The representation of $R(\la)$ in \eqref{eq_rationalmatrix} and the block Kronecker linearizations $S(\la)$ of $R(\la)$ \eqref{block_kronecker_lin} are the two fundamental ingredients of this paper.

As explained in \cite[Section 3.1]{ADMZ}, the finite eigenvalues, together with their partial multiplicities, of $S(\la)$ (resp. $A-\la I_\ell$) coincide with the finite zeros (resp. poles) of $R(\la)$, together with their partial multiplicities. Moreover, the eigenvalue structure at infinity of $S(\la)$ allows us to obtain via a simple shift rule the pole-zero structure at infinity of $R(\la)$.\footnote{More precisely, according to \cite[p. 1683]{ADMZ}, if $r$ is the normal rank of $R(\la)$ and $e_1 \leq \cdots \leq e_r$ are the $r$ largest partial multiplicities at infinity of $S(\la)$, then $e_1 -d \leq \cdots \leq e_r-d$ are the structural indices at infinity of $R(\la)$.} In addition, as proved in \cite[Section 6]{ADMZ-2}, the right (resp. left) minimal indices of $S(\la)$ are those of $R(\la)$ plus $\epsilon$ (resp. $\eta$). Thus, $S(\la)$ comprises the complete eigenstructure of $R(\la)$. Observe that the application to $S(\la)$ of the QZ algorithm \cite{MOS}, in the regular case, or of the staircase algorithm \cite{VanDooren79}, in the singular case, gives the zeros and the minimal indices, in the singular case, of $R(\la)$, but not the poles, which are in $A-\la I_\ell$.

It is worth mentioning that although the families of block Kronecker linearizations of polynomial \cite{DLPV} and rational \cite{ADMZ} matrices are very recent, some particular examples of strong linearizations in these families appeared much earlier in the literature. For instance, it was shown in \cite{VDD83} that a valid ``realization'' for the polynomial part $D(\la)$ in \eqref{eq_rationalmatrix} is given by the following minimal Rosenbrock polynomial system matrix \cite{Ros70}
$$ S_D(\la) : =
\left[ \begin{array}{cccc|c} I_n & -\la I_n & & & \\  & I_n & \ddots & &   \\ & & \ddots & -\la I_n \\ & & & I_n & -\la I_n \\
\hline \la D_d & \ldots & \ldots & \la D_2 & \la D_1+D_0  \end{array}\right]:= \left[ \begin{array}{c|c} T(\la) & -U(\la) \\ \hline V(\la) & W(\la) \end{array}\right] ,
$$
which means that $D(\la)=W(\la)+V(\la)T(\la)^{-1}U(\la)$. It is easy to see that after moving the bottom block row of $S_D (\la)$ to the top position, a block Kronecker linearizarion of $D(\la)$ is obtained with $K_2 (\la)$ empty \cite[Section 4]{DLPV}. Combining the minimal state-space realization $C(\la I_\ell-A)^{-1}B$ and the polynomial system matrix $S_D(\la)$ yields the following minimal polynomial system matrix for the rational matrix $R(\la)$ in \eqref{eq_rationalmatrix} :
{\small$$ S_R(\la) := \left[ \begin{array}{c:cccc|c} A-\la I_\ell & & & & &  B \\ \hdashline
& I_n & -\la I_n & & & \\ & & I_n & \ddots & & \\ & & & \ddots & -\la I_n \\ & & & & I_n & -\la I_n \\
\hline C & \la D_d & \ldots & \ldots & \la D_2 & \la D_1+D_0  \end{array}\right] := \left[ \begin{array}{c|c} T_R(\la) & -U_R(\la) \\ \hline V_R(\la) & W_R(\la) \end{array}\right] ,
$$}
i.e., $R(\la)=W_R(\la)+V_R(\la)T_R(\la)^{-1}U_R(\la)$. A pencil with a structure similar to $ S_R(\la)$ can also be found in \cite{SuBai}. It is easy to see that, modulo some block permutations, $S_R(\la)$ is a particular case of the block Kronecker linearizations appearing in \eqref{block_kronecker_lin} for $R(\la)$, with $K_2 (\la)$ empty and $\widehat{K}_2 = I_m$. In fact, it can be proved that what has been shown above for $S_R (\la)$ holds for all the block Kronecker linearizations of $R(\la)$ in \eqref{block_kronecker_lin}, since all of them can be seen as minimal Rosenbrock polynomial system matrices of $R(\la)$ when we permute them to
$$ S_K(\la) := \left[ \begin{array}{c:ccc} A-\la I_\ell & 0 &  B\widehat K_1  \\ \hdashline
0 & 0 & K_1(\la) \\ \widehat K_2^TC & K_2^T(\la) & M(\la) \end{array}\right] $$
and, then, we partition them appropriately, since the bottom right submatrix is a linearization for $D(\la)$. This approach based on the
block Kronecker linearizations for polynomial matrices, also contains the companion forms as a special case.

It was shown in \cite{DLPV} that perturbations of the block Kronecker linearizations of a polynomial matrix $D(\la)$ can be mapped to perturbations of the coefficients of $D(\la)$ without significant growth of the relative norms of the perturbations under mild assumptions that require to scale $D(\la)$ to have norm equal to $1$ and to use linearizations with the norm of $M(\la)$ of the same order as the norm of $D(\la)$ (see \cite[Corollary 5.24]{DLPV}). As a corollary of this perturbation result, we obtain that under such assumptions the computation of the eigenvalues and minimal indices of a polynomial matrix by applying the $QZ$ or the staircase algorithm to one of its block Kronecker linearizations is a backward stable method from the point of view of the polynomial matrix. In this paper we show that this can be extended to rational matrices as well, considering as coefficients of the rational matrix those in the quadruple $\{\la I_\ell-A,B,C,D(\la)\}$. However, we emphasize that the perturbation analysis for block Kronecker linearizations of rational matrices is considerably more complicated than the one in \cite{DLPV} and, therefore, we limit ourselves to perform a first order analysis. We also remark that the scaling needed to get satisfactory perturbation bounds is more delicate than the one in \cite{DLPV}. As far as we know, this is the first structural backward error analysis of this type performed in the literature for linearizations of rational matrices.

We assume throughout the paper that $\ell > 0$ since, otherwise, $R(\la)$ in \eqref{eq_rationalmatrix} is a polynomial matrix and this case was studied in \cite{DLPV}. Except in Subsection \ref{subsec.linearD}, we also assume that at least one of the parameters $\epsilon$ and $\eta$ in \eqref{block_kronecker_lin} is larger than zero since, otherwise, none of the blocks $K_1 (\la)$ and $K_2(\la)$ appears and block Kronecker linearizations collapse to much simpler pencils. Note that $\max (\eta,\epsilon) >0$ implies that the degree $d$ of the polynomial part $D(\la)$ of $R(\la)$ is larger than $1$. The simple case $d \leq 1$ is studied in Subsection \ref{subsec.linearD}.

In order to measure perturbations, we need to introduce appropriate norms for pencils, polynomial matrices and rational matrices expressed as in \eqref{eq_rationalmatrix}. For any pair of matrices $X$ and $Y$
of arbitrary dimensions (that might be different), we will use the following norms
$$ \|(X,Y)\|_F := \left(\|X\|_F^2 +\|Y\|^2_F\right)^{\frac12} = \|\, [\mathrm{vec}(X)^T,\mathrm{vec}(Y)^T ] \, \|_2,$$
$$ \|(X,Y)\|_2 := \left(\|X\|_2^2 +\|Y\|^2_2\right)^{\frac12},$$
where $\|X\|_F$ and $\|X\|_2$ are, respectively, the Frobenius and spectral matrix norms and $\mathrm{vec}(X)$ is the operator that stacks the columns of a matrix into one column vector \cite{higham-book-2002}. For a pencil $S(\la):=A-\la B$  we define the corresponding norms via the two matrix coefficients~:
$$  \|S(\la)\|_F := \| (A,B) \|_F, \quad \|S(\la)\|_2 := \| (A,B) \|_2.$$
More generally, for a polynomial matrix $D(\la):=\sum_{i=0}^d D_i\la^i$, we will use the norm
$$ \|D(\la)\|_F := \sqrt{ \sum_{i=0}^d \| D_i\|_F^2},$$
and for a list of polynomial matrices $(D_{1}(\la),\ldots,D_{p}(\la))$, the norm
$$ \|(D_{1}(\la),\ldots,D_{p}(\la))\|_F := \sqrt{ \sum_{i=1}^p \| D_{i}(\la)\|_F^2}.$$ Finally, for a rational matrix $R(\la),$ represented by a quadruple $\{\la I_\ell- A,B,C,D(\la)\},$ as in \eqref{eq_rationalmatrix}, we use
the ``norm'' $$ \|R(\la)\|_F :=\|(\la I_\ell-A,B,C,D(\la))\|_F= \sqrt{\ell + \| A \|_F^2 + \| B \|_F^2+ \| C \|_F^2+ \sum_{i=0}^d\| {D_i}\|_F^2}.$$
That is, the ``norm'' of a rational matrix $R(\la)$ is defined as the norm of an associated polynomial system matrix $P(\la)$, in this case,
\begin{equation}\label{def:norm}
\|R(\la)\|_F := \|P(\la)\|_F \quad \text{where}\quad P(\la):=\begin{bmatrix}
\la I_\ell-A & -B\\
C & D(\la)
\end{bmatrix}.
\end{equation}
We remark that $\|R(\la)\|_F$ is not rigorously a ``norm'' for $R(\la)$ because, for instance, $R(\la)$ is zero if $B=0$ and $D(\la) = 0$, but $\|R(\la)\|_F$ is not. Despite this fact, and with a clear abuse of nomenclature, we will use the terminology ``norm of a rational matrix'' in the sense explained above.

The paper is organized as follows. After this introductory section, we describe in Section \ref{Sylvester} the basic systems of matrix equations we will use
in this paper, and, in Section \ref{SVbounds}, some bounds for the singular values of certain matrices related to these systems of matrix equations. In Section \ref{structure} we explain how to restore the structure of block Kronecker  linearizations of rational matrices after they suffer sufficiently small perturbations, and, in Section \ref{scaling}, we derive a scaling technique that allows us to guarantee structured backward stability for (regular or singular) rational eigenvalue problems solved via block Kronecker linearizations.
Finally, in Section \ref{numerical} we give a number of numerical results illustrating our theoretical bounds and, in Section \ref{sec.conclusions}, we establish some conclusions.

\section{Generalized Sylvester equations} \label{Sylvester}

In order to restore the structure of perturbed block Kronecker linearizations of rational matrices, we will need to guarantee that some matrix equations have solutions and to bound the norm of their minimal norm solution. The matrix equations that we will encounter are particular cases of the generalized Sylvester equation for  $m_i\times n_i$ pencils of matrices  $A_i-\la B_i$, $i=1,2$, which is the following equation
in the unknowns $X$ and $Y$~:
\begin{equation} \label{GS}  X(A_1-\lambda B_1)+(A_2-\la B_2)Y = \Delta^a-\la \Delta^b.
\end{equation}
It is easily seen to be equivalent to a linear system of equations, when rewriting it as
\begin{eqnarray*}
	XA_1+A_2Y & = & \Delta^a,\\ XB_1+B_2Y & = & \Delta^b,
\end{eqnarray*}
or, when using Kronecker products and the vec$(\cdot)$ notation, as
\begin{equation} \label{vec}
\left[
\begin{array}{c|c}
A_1^T\otimes I_{m_2} & I_{n_1} \otimes A_2 \\\hline
B_1^T\otimes I_{m_2} & I_{n_1} \otimes B_2
\end{array}
\right]
\left[
\begin{array}{r}
\mathrm{vec}(X)\\\hline
\mathrm{vec}(Y)
\end{array}
\right]= \left[
\begin{array}{c}
\mathrm{vec}(\Delta^a)\\\hline
\mathrm{vec}(\Delta^b)
\end{array}
\right].
\end{equation}
The dimension of the unknowns $X$ and $Y$ are $m_2\times m_1$ and $n_2\times n_1$, respectively,
and those of the right hand sides $\Delta^a$ and $\Delta^b$ are each $m_2\times n_1$. These equations will be used in this paper in two contexts, which we briefly recall here.

\medskip

{\bf Block elimination.}
Let $A_i-\la B_i$ be two $m_i\times n_i$ pencils, $i=1, 2$, that have respectively full column normal rank $n_1$ and full row normal rank $m_2$.
Then the problem of block anti-diagonalizing  the pencil {\small $ \left[\begin{array}{cc} 0 &  A_1-\lambda B_1 \\  A_2-\la B_2  & \Delta^a- \la \Delta^b \end{array}\right]$}, that is, finding $X$ and $Y$ such that
{\small \begin{equation} \label{GS2}
	\left[\begin{array}{cc} I_{m_1} & 0 \\ -X & I_{m_2} \end{array} \right]
	\left[\begin{array}{cc} 0 &  A_1-\lambda B_1 \\  A_2-\la B_2  & \Delta^a- \la \Delta^b \end{array}\right]
	\left[\begin{array}{cc} I_{n_2} & -Y  \\ 0 & I_{n_1} \end{array} \right] =
	\left[\begin{array}{cc} 0 &  A_1-\lambda B_1 \\ A_2-\la B_2 & 0 \end{array}\right],
	\end{equation}}
amounts to finding a solution for the generalized Sylvester equation \eqref{GS}.
It is known that there exists a solution $(X,Y) \in \mathbb{F}^{m_2\times m_1} \times \mathbb{F}^{n_2\times n_1}$  for a {\em particular} right hand side $(\Delta^a,\Delta^b) \in \mathbb{F}^{m_2\times n_1} \times \mathbb{F}^{m_2\times n_1}$ {\em if and only if} the pencils
$$ \left[\begin{array}{cc} 0 &  A_1-\lambda B_1 \\  A_2-\la B_2  & \Delta^a- \la \Delta^b \end{array}\right] \quad \mathrm{and} \quad
\left[\begin{array}{cc} 0 &  A_1-\lambda B_1 \\ A_2-\la B_2 & 0 \end{array}\right]$$
are strictly equivalent (i.e. have the same Kronecker structure) \cite{Andrej}. But in order to have a solution {\em for any right hand side}
$\Delta^a-\la \Delta^b$ one requires the stronger condition that the pencils $A_1-\lambda B_1$ and $A_2-\lambda B_2$ have no common generalized eigenvalues (see \cite{VanDooren83}). We recall here the result proven in \cite{VanDooren83} that is relevant for our work.

\begin{theorem}[\cite{VanDooren83}]
	\label{th:GS2} Let the pencils  $A_i-\la B_i$ of dimensions $m_i\times n_i, i=1,2$, be respectively of  full column normal rank $n_1\le m_1$ and of
	full row normal rank $m_2\le n_2$, and let these two pencils have no common generalized eigenvalues. Then there always exists a solution $(X,Y)$ to the system of equations \eqref{GS2}, for any perturbation $\Delta^a- \la \Delta^b$. Moreover, the generalized eigenvalues of the
	pencil \eqref{GS2} are the union of the generalized eigenvalues of the pencils  $A_i-\la B_i, \; i=1,2.$
\end{theorem}

The system is underdetermined if either of the two inequalities  $m_1 \ge n_1$ and $n_2\ge m_2$, is strict. Under the hypotheses of Theorem \ref{th:GS2}, the system \eqref{vec} must be compatible for any right hand side, and hence the Kronecker product matrix in the left hand side of \eqref{vec} must have full row rank $2m_2n_1$.
A bound for the minimum Frobenius-norm solution $(X,Y)$ is then obtained in terms of the smallest singular value $\sigma_{2m_2n_1}$ of the matrix in \eqref{vec}~:
\begin{equation} \label{GSbound} \|(X,Y)\|_F \le \frac{\|(\Delta^a,\Delta^b)\|_F}{ \sigma_{2m_2n_1}\left(	\left[
	\begin{array}{c|c}
	A^T_1\otimes I_{m_2} &I_{n_1} \otimes A_2 \\\hline
	B^T_1\otimes I_{m_2} & I_{n_1} \otimes B_2
	\end{array}
	\right]\right) }.
\end{equation}

\medskip

{\bf Equivalent pencils.} The second problem in this paper where a generalized Sylvester equation as in \eqref{GS} arises is that of strictly equivalent pencils (see e.g. \cite{Gan}). Let the pencils  $A_i-\la B_i, \; i=1,2$, be both of dimension $m\times n$, then they are strictly equivalent if and only if there exist invertible matrices $S$ and $T$ such that $S(A_1-\la B_1)=(A_2-\la B_2)T$. Such pencils must then have the same Kronecker canonical form \cite{Gan}. We are interested in finding the solution where $S$ and $T$ are as close as possible to the identity matrix.
This can be achieved by writing the transformation matrices as
$$  S=I+X, \quad T =I-Y $$
and then minimizing the Frobenius norm of the pair $(X,Y)$. The corresponding equations are then $$ (I+X)(A_1-\la B_1)=(A_2-\la B_2)(I-Y)$$
or, when putting $\Delta^a-\la \Delta^b := (A_2-\la B_2)-(A_1-\la B_1)$, we finally obtain
\begin{equation}\label{GS1}
X(A_1-\la B_1)+(A_2-\la B_2)Y = \Delta^a-\la \Delta^b,
\end{equation}
which is again solved by using \eqref{vec}.
We will use this to ``restore'' a slightly perturbed pencil $(A_2-\la B_2):=(A_1-\la B_1) + (\Delta^a-\la \Delta^b)$ to its original form $(A_1-\la B_1)$
using a strict equivalence transformation
\begin{equation} \label{GS1aa} (I+X)^{-1}(A_2-\la B_2)(I-Y)=A_1-\la B_1
\end{equation}
that is very close to the identity, when we are sure that both pencils have the same Kronecker canonical form. The bounds for the norm of $X$ and $Y$ are in fact given by \eqref{GSbound} for which we derive exact expressions in the next section.
Notice that we can not apply Theorem \ref{th:GS2} to prove existence of a solution for equation \eqref{GS1},  since in this case both pencils must have the same generalized eigenvalues and the same normal rank. A sufficient condition for the consistency of \eqref{GS1} is that $A_1-\la B_1$ and $A_2-\la B_2$ have the same Kronecker canonical form.

The condition that the Kronecker canonical form of a pencil does not change under arbitrary sufficiently small perturbations only holds for very special pencils. In particular, it holds for the Kronecker product of Kronecker blocks times identity matrices, i.e., for $L_k(\la)\otimes I_r$. This is a consequence of the results in \cite{VanD17}, because $L_k(\la)\otimes I_r$ has full-Sylvester-rank by \cite[Theorem 4.3(a)]{VanD17} and, then, \cite[Theorem 6.6]{VanD17} guarantees that $L_k(\la)\otimes I_r + (\Delta^a-\la \Delta^b)$  has the same Kronecker canonical form as $L_k(\la)\otimes I_r$ for all the perturbations $(\Delta^a,\Delta^b)$ whose norms are smaller than the bounds in \cite[Theorem 6.6]{VanD17}. Since we will solve \eqref{GS1}-\eqref{GS1aa} only in the case $A_1 - \la B_1 = L_k(\la)\otimes I_r$, these results prove that
\eqref{GS1} has a solution for all sufficiently small perturbations $(\Delta^a,\Delta^b)$ in the cases of interest in this paper.

\section{Singular value bounds} \label{SVbounds}

In the analysis of Section \ref{structure}, we will need upper bounds for the minimum norm solutions of the generalized Sylvester equation \eqref{GS} for pairs of pencils $(A_i-\la B_i), \; i=1,2$, which all involve Kronecker blocks
$L_k(\la):=E_k-\la F_k$, where the $k\times (k+1)$ matrices $E_k$ and $F_k$ are given by
$$ E_k:= \left[ \begin{array}{cccccc}  1 & 0 & \\ &  1 & 0 \\
& & \ddots & \ddots \\ & & & 1 & 0 \end{array} \right]  \quad \mathrm{and} \quad
F_k:= \left[ \begin{array}{cccccc}  0 & 1 & \\ &  0 & 1 \\
& & \ddots & \ddots \\ & & &  0 & 1 \end{array} \right].
$$
To find such upper bounds is equivalent to find lower bounds for the singular values in the denominator of the right hand side of \eqref{GSbound}.
We consider the generalized Sylvester equations for the following list of pencil pairs with their smallest singular value of the corresponding linear maps:

\begin{enumerate}
	\item $A_1 - \la B_1 =  A-\la I_\ell$ and $A_2 - \la B_2 =  L_\e(\la)\otimes I_n$:
	\begin{equation}\label{singval_1}
	\omega_1 :=
	\sigma_{2\ell \e n}\left[ \begin{array}{c|c} A^T\otimes I_{\e n} & I_{\ell} \otimes E_\e\otimes I_n \\ \hline
	I_\ell \otimes I_{\e n} & I_{\ell} \otimes F_\e\otimes I_n \end{array} \right].
	\end{equation}
	
	\item $A_1 - \la B_1 =  L_\eta^T(\la)\otimes I_m$ and $A_2 - \la B_2 = A-\la I_\ell$:
	\begin{equation}\label{singval_2}
	\omega_2 :=
	\sigma_{2\eta m \ell}\left[ \begin{array}{c|c} E_\eta \otimes I_{m \ell} & I_{\eta m} \otimes A \\ \hline
	F_\eta \otimes I_{m\ell} & I_{\eta m} \otimes I_\ell \end{array} \right].
	\end{equation}
	
	\item $A_1 - \la B_1 =  L_\eta^T(\la)\otimes I_m$ and $A_2 - \la B_2 = L_\e(\la)\otimes I_n$:
	\begin{equation}\label{singval_3}
	\omega_3 :=
	\sigma_{2\eta m\e n}\left[ \begin{array}{c|c} E_\eta\otimes I_{m\e n} & I_{\eta m} \otimes E_\e\otimes I_n \\ \hline
	F_\eta \otimes I_{m \e n} & I_{\eta m} \otimes F_\e\otimes I_n \end{array} \right] .
	\end{equation}
	\item $A_1 - \la B_1 =  L_k(\la)\otimes I_r$ and $A_2 - \la B_2 =  L_k(\la)\otimes I_r$:
	\begin{equation}\label{singval_4}
	\omega_4 :=
	\sigma_{2(k+1)rkr}\left[ \begin{array}{c|c} E_k^T \otimes I_{rkr} & I_{(k+1)r} \otimes E_k\otimes I_r \\ \hline
	F_k^T \otimes I_{rkr} & I_{(k+1)r} \otimes F_k\otimes I_r \end{array} \right].
	\end{equation}
\end{enumerate}

In Lemma \ref{lem:singval_1} we analyze the first problem and give a lower bound for $\omega_1$.

\begin{lemma}\label{lem:singval_1} Let $\omega_1$ be the singular value in \eqref{singval_1}. Then
	
	\begin{equation}\label{eq:omega1}
	\omega_1 \ge \frac{1}{1+2\e\max(1,\|A\|_2^\e)}.
	\end{equation}
	
\end{lemma}
\begin{proof}
	It follows from the properties of singular values of Kronecker products that $\omega_1$ is also equal to
	$$ \omega_1 =
	\sigma_{2\ell \e}\left[ \begin{array}{c|c}
	A^T\otimes I_{\e} & I_{\ell} \otimes E_\e \\ \hline
	I_\ell \otimes I_{\e} & I_{\ell} \otimes F_\e
	\end{array} \right]
	$$
	and using perfect shuffle permutations we also get
	$$  \omega_1=
	\sigma_{2\e \ell}\left[ \begin{array}{c|c}
	I_{\e} \otimes A^T & E_\e \otimes I_\ell \\ \hline
	I_\e \otimes I_{\ell} & F_{\e} \otimes I_\ell
	\end{array} \right].
	$$
	The smallest singular value $ \sigma_{2\e\ell}$ is larger than the smallest singular value of any $2\e \ell \times 2\e \ell$ submatrix.
	Let us take for this the submatrix obtained by dropping the last block column~:
	$$ M =\left[ \begin{array}{c|c}
	I_{\e} \otimes A^T & I_\e \otimes I_\ell \\ \hline
	I_\e \otimes I_{\ell} & J_{\e} \otimes I_\ell
	\end{array} \right],  \quad \mathrm{where} \quad J_\e := \left[ \begin{array}{ccccc}  0 & 1 & \\ &  0 & \ddots \\
	& & \ddots & 1 \\ & & &  0  \end{array} \right]  \in \mathbb{F}^{\e \times \e} .$$
	We can factorize this matrix as
	$$ M =\left[ \begin{array}{c|c}
	I_{\e} \otimes A^T & I_\e \otimes I_\ell \\ \hline
	I_\e \otimes I_{\ell} & 0
	\end{array} \right]\left[ \begin{array}{c|c}
	I_{\e \ell} & 0  \\ \hline
	0 &  I_{\e \ell}- J_{\e} \otimes A^T
	\end{array} \right]\left[ \begin{array}{c|c}
	I_{\e} \otimes I_\ell & J_{\e} \otimes I_\ell \\ \hline
	0  & I_\e\otimes I_\ell
	\end{array} \right].
	$$
	Therefore its inverse equals
	\begin{align*} M^{-1} & =\left[ \begin{array}{c|c}
	I_{\e} \otimes I_\ell & - J_{\e} \otimes I_\ell \\ \hline
	0  & I_\e\otimes I_\ell
	\end{array} \right]\left[ \begin{array}{c|c}
	I_{\e \ell} & 0  \\ \hline
	0 &  (I_{\e \ell}- J_{\e} \otimes A^T)^{-1}
	\end{array} \right]
	\left[ \begin{array}{c|c}
	0 & I_\e \otimes I_\ell \\ \hline
	I_\e \otimes I_{\ell} &  -I_{\e} \otimes A^T
	\end{array} \right] \\
	& =\left[ \begin{array}{c}
	I_{\e \ell}  \\ \hline 0
	\end{array} \right]
	\left[ \begin{array}{c|c}
	0 & I_{\e \ell}
	\end{array} \right] +\left[ \begin{array}{c}
	- J_{\e} \otimes I_\ell \\ \hline
	I_{\e \ell}
	\end{array} \right]  (I_{\e \ell}- J_{\e} \otimes A^T)^{-1}
	\left[ \begin{array}{c|c}
	I_{\e \ell} &  -I_{\e} \otimes A^T
	\end{array} \right].
	\end{align*}
	It then follows that
	$$ \|M^{-1}\|_2\le 1 + \sqrt{2}\sqrt{1+\|A\|^2_2}\left[1 + \|A\|_2 + \|A\|_2^2 + \ldots + \|A\|_2^{\e-1}\right],$$
	since
	$$
	(I_{\e \ell}- J_{\e} \otimes A^T)^{-1} = \sum_{i=0}^{\e-1}  J_{\e}^i \otimes {A^i}^T.
	$$
	In particular, for $\|A\|_2\le 1$ we obtain the bound $\|M^{-1}\|_2 \le 1 + 2\e$, while  for  $\|A\|_2 > 1$ we obtain the bound
	$\|M^{-1}\|_2 \le 1 + 2\e\|A\|_2^\e$. This finally yields the inequality
	\begin{equation*}
	\omega_1 \ge \frac{1}{1+2\e\max(1,\|A\|_2^\e)}.
	\end{equation*}
\end{proof}

The second generalized Sylvester equation is essentially the transposed of the first equation and the analysis is therefore completely analogous. This immediately yields Lemma \ref{lem:singval_2}.
\begin{lemma}\label{lem:singval_2} Let $\omega_2$ be the singular value in \eqref{singval_2}. Then
	\begin{equation}\label{eq:omega2}
	\omega_2 \ge \frac{1}{1+2\eta\max(1,\|A\|_2^\eta)}.
	\end{equation}
\end{lemma}

The third generalized Sylvester equation was analyzed in \cite{DLPV} and its associated smallest singular value is exactly equal to
$\omega_3 = 2\sin(\pi/(4\min(\e,\eta)+2))$ if $\e\neq \eta$, and to $ 2\sin(\pi/4\eta)$ if $\e=\eta$.
Notice that we can assume $\min(\e,\eta)\ge 1$ since otherwise the equation is void.
For $\e\neq \eta$ we then obtain $\omega_3\ge \frac{3}{2\min(\e,\eta)+1}$ since $\sin x\ge 3x/\pi$ for $0 \leq x\le \pi/6$, and
for $\e=\eta$ we then obtain $\omega_3\ge\frac{\sqrt{2}}{\eta} $ since $\sin x\ge 2\sqrt{2}x/\pi$ for $0 \leq x\le \pi/4$.
We have also that $2\eta= \e+\eta$ if $\e=\eta$ and $2\min(\e,\eta)+1\le \e+\eta$ if $\e\neq \eta$, which finally yields the lower bound in Lemma \ref{lem:singval_3} for $\omega_3.$
\begin{lemma}\label{lem:singval_3} Let $\omega_3$ be the singular value in \eqref{singval_3}. Then
	\begin{equation}\label{eq:omega3}
	\omega_3 \ge \frac{2\sqrt{2}}{(\e+\eta)}.
	\end{equation}
\end{lemma}

In Lemma \ref{lem:singval_4}, we give a lower bound for the smallest singular value $\omega_4$ corresponding to the fourth generalized Sylvester equation.

\begin{lemma}\label{lem:singval_4} Let $\omega_4$ be the singular value in \eqref{singval_4}. Then
	\begin{equation}\label{eq:omega4}
	\omega_4 \ge \frac{3}{4k-1}.
	\end{equation}
\end{lemma}
\begin{proof}
	We prove first that $\omega_4=2\sin(\pi/(8k-2))$.
	This is obtained as follows. We can again use the properties of Kronecker products to prove that
	$$ \omega_4 =  \sigma_{2k(k+1)}\left[ \begin{array}{c|c} E_k^T \otimes I_{k} & I_{(k+1)} \otimes E_k \\ \hline
	F_k^T \otimes I_{k} & I_{(k+1)} \otimes F_k \end{array} \right].
	$$
	This matrix can be transformed by row and column permutations to the direct sum of smaller matrices~:
	$$  M_1\oplus M_1\oplus M_3 \oplus M_3 \oplus \cdots \oplus M_{2k-1}  \oplus  M_{2k-1} \oplus N_{2k} ,
	$$
	see \ref{Appendix}, where the blocks
	\begin{equation} \label{MkNk}
	M_k := \left[ \begin{array}{ccccc}  1 & 1 & \\ &  1 & \ddots \\
	& & \ddots & 1 \\ & & &  1  \end{array} \right] \in \mathbb{F}^{k\times k} , \quad N_k  := \left[ \begin{array}{cccccc}  1 & 1 & \\ &  1 & \ddots \\
	& & \ddots & 1 \\ & & & 1 & 1 \end{array} \right]  \in \mathbb{F}^{k\times (k+1)}
	\end{equation}
	have as smallest singular values $2\sin \frac{\pi}{4k+2}$ and  $2\sin \frac{\pi}{2k+2}$, respectively (see \cite[Proof of Proposition B.4]{DLPV}). The smallest singular value therefore corresponds to
	$M_{2k-1}$ and equals $\omega_4=2\sin(\pi/(8k-2))$. For $k \ge 1$, we use again that $\sin x \ge 3x/\pi$ for $0 \leq x\le \pi/6$, to obtain the bound
	$\omega_4 \ge \frac{3}{4k-1}.
	$\end{proof}

\section{Restoring the rational structure of the linearization after perturbations} \label{structure}

We now consider perturbations of the following block Kronecker linearization introduced in \eqref{block_kronecker_lin}
\begin{equation}\label{eq:pencil}
S(\la):=\left[\begin{array}{ccc} S_{11}(\la) & S_{12}(\la)  & S_{13}(\la) \\ S_{21}(\la)  & S_{22}(\la)  & 0\\ S_{31}(\la) &  0 &  0 \end{array}\right] := \left[\begin{array}{ccc} M(\la) & \widehat K_2^T C &  K_2^T(\la)  \\ B \widehat K_1 & A- \la I_\ell & 0\\
K_1(\la) &   0 &  0
\end{array}\right],
\end{equation}
where $S_{13}(\la)$ is $(\eta+1)m \times \eta m$ and  has full column rank $\eta m$,  $S_{22}(\la)$ is $\ell \times \ell$ and is a regular pencil, $S_{31}(\la)$ is $\e n \times (\e+1)n$ and has full row rank $\e n$, and where no two of these three pencils have common generalized eigenvalues.  As explained in the introduction, if the state-space triple $\{A,B,C\}$ is minimal, then $S(\la)$ is a strong linearization of the $m\times n$ rational matrix
\begin{equation} \label{eq.Rsect4fro}
R(\la)= C(\la I_\ell-A)^{-1}B + (\Lambda_\eta(\la)\otimes I_m)^T M(\la) (\Lambda_\epsilon(\la)\otimes I_n).
\end{equation}
Except in Subsection \ref{subsec.linearD}, we assume in this section that $\max (\eta,\epsilon) > 0$. This means that the degree $d = \epsilon + \eta + 1$ of the polynomial part $D(\la) = (\Lambda_\eta(\la)\otimes I_m)^T M(\la) (\Lambda_\epsilon(\la)\otimes I_n)$ of $R(\la)$ is greater than $1$ and that at least one of the blocks $K_1(\la)$ or $K_2(\la)$ is not an empty matrix. The degenerate case in which $\e=0$ and $\eta=0$ will be studied in Subsection \ref{subsec.linearD}.

Since $S(\la)$ is a strong linearization of $R(\la)$, $S(\la)$ has the exact eigenstructure of the finite zeros of $R(\la)$, and its infinite zero structure as well as its left and right null-space structure can be correctly retrieved from the pencil via simple constant shifts, as explained in the introduction.  In order to  compute this eigenstructure, we make use of the staircase algorithm \cite{VanDooren79}, followed by the $QZ$ algorithm \cite{MOS}, on $S(\la)$. The backward stability of these two algorithms guarantees in fact that we computed the exact eigenstructure of a slightly perturbed pencil
\begin{equation}\label{eq:perturbed}
\widehat S(\la):= S(\la) + \Delta_S(\la), \quad  \Delta_S(\la) := \left[\begin{array}{ccc} \Delta_{11}(\la) & \Delta_{12}(\la) &  \Delta_{13}(\la)  \\ \Delta_{21}(\la) & \Delta_{22}(\la) & \Delta_{23}(\la)\\
\Delta_{31}(\la) &  \Delta_{32}(\la) &  \Delta_{33}(\la)
\end{array}\right],
\end{equation}
where the pencil $\Delta_S(\la)$ has a norm which is much smaller than the norm of $S(\la)$. More precisely, $\|\Delta_S(\la)\|_F = O(\epsilon_M) \, \|S(\la)\|_F$, where $\epsilon_M$ is the machine precision of the computer. But even for very small perturbations, the structure
of the pencil $S(\la)$ is lost, and therefore also the connection between $\widehat S(\la)$ and some rational matrix $\widehat R(\la)$ is lost. In this section, we will show that this structure can be restored, without affecting the computed eigenstructure.  For this, one needs only to find a strict equivalence transformation that is close to the identity and restores the structure of $\wh{S}(\la)$ to a new pencil $\widetilde{S}(\la)$ that is a block Kronecker linearization, with the same parameters $\epsilon$ and $\eta$ as $S(\la)$, of a rational matrix $\widetilde R(\la)$~:
\begin{equation} \label{eqstep33}
\widetilde S(\la):=(I-X)(S(\la)+\Delta_S(\la))(I-Y) = \left[\begin{array}{ccc} \widetilde M(\la) & \widehat K_2^T \widetilde C &  K_2^T(\la)  \\ \widetilde B \widehat K_1 & \widetilde A- \la I_\ell & 0\\
K_1(\la) &   0 &  0
\end{array}\right].
\end{equation}
We will see that if $\|\Delta_S(\la)\|_F$ is sufficiently small, then the perturbed system triple $\{ \widetilde A,\widetilde B,\widetilde C\}$ is very close to the unperturbed minimal one $\{ A, B, C\}$ and, so, $\{ \widetilde A,\widetilde B,\widetilde C\}$  is still minimal, since minimality is a generic property equivalent to the controllability matrix having full row rank and the  observability matrix having full column rank \cite[Chapter 6]{Kai80}. Observe that according to \cite{ADMZ}, or the discussion in the introduction, $\widetilde{S} (\la)$ is a strong linearization of the $m \times n$ rational matrix
\begin{equation}\label{nearby_rational}
\begin{split}
\widetilde{R} (\la) &:= 	\widetilde C(\la I_\ell-	\widetilde A)^{-1}	\widetilde B + (\Lambda_\eta(\la)\otimes I_m)^T \widetilde M(\la) (\Lambda_\epsilon(\la)\otimes I_n)  \\ & =: \widetilde C(\la I_\ell-	\widetilde A)^{-1}	\widetilde B + \widetilde D (\la) \, .
\end{split}
\end{equation}
Since the eigenstructures of the pencils $\widehat{S} (\la)$ and $\widetilde{S} (\la)$ are identical, the results in this section prove that the computed finite eigenvalues of $S(\la)$ and their partial multiplicities are the exact finite zeros and their partial multiplicities of $\widetilde{R} (\la)$, the computed right (resp. left) minimal indices of $S(\la)$ minus $\epsilon$ (resp. $\eta$) are the exact right (resp. left) minimal indices of $\widetilde{R} (\la)$, and, if a number $\nu_r$ of right minimal indices of $S(\la)$ have been computed, then the computed $n-\nu_r$ largest partial multiplicities at infinity of $S(\la)$ minus $d$ are the exact structural indices at infinity of $\widetilde{R} (\la)$. This is a very strong backward error result for the computation of the eigenstructure of $R(\la)$ in the case we are able to prove that $\| \widetilde{A} - A\|_F,  \|\widetilde{B} - B\|_F, \| \widetilde{C} - C\|_F$ and $\| \widetilde{D} (\la) - D(\la) \|_F$ are very small.

The restoration of the structure in $\widehat{S} (\la)$ will be done in three steps, each of them involving a strict equivalence transformation close to the identity:
\begin{itemize}
	\item  \textbf{Step 1:} We restore the block anti-triangular structure of the perturbed pencil $\wh{S}(\la),$ i.e., the blocks (2,3), (3,2) and (3,3) are transformed to become 0.
	\item  \textbf{Step 2:} We take care of the anti-diagonal blocks (1,3), (2,2) and (3,1), by restoring their $0$ and $I$ block matrices.
	\item \textbf{Step 3:} We restore the special structure of the blocks (1,2) and (2,1).
\end{itemize}
At each step $k,$ for $k=1,2,3,$ we obtain a pencil
\begin{equation} \label{eq.stepkstrictequiv}
\wh{S}_{k}(\la):=(I-X_k)\wh{S}_{k-1}(\la)(I-Y_k):=\wh{S}_{k-1}(\la)+\Delta_{k}(\la),
\end{equation}
where $\wh{S}_{0}(\la):=\wh S(\la)$ and $\Delta_0(\la):=\Delta_S(\la)$~:
$$ S(\la) \stackrel{+\Delta_{0}(\la)}{\xrightarrow{\hspace{0.6cm}}} \wh S(\la)=\wh S_0(\la) \stackrel{+\Delta_{1}(\la)}{\xrightarrow{\hspace{0.6cm}}} \wh S_1(\la) \stackrel{+\Delta_{2}(\la)}{\xrightarrow{\hspace{0.6cm}}} \wh S_2(\la) \stackrel{+\Delta_{3}(\la)}{\xrightarrow{\hspace{0.6cm}}}\wh S_3(\la)= \widetilde S(\la).
$$

We will compute bounds for $\|(X_k,Y_k)\|_F$ as a function of $\|\wh S_{k-1}(\la)\|_F$,
where the Frobenius norms are computed as defined in the introduction. Moreover, we define the cumulative errors:
\begin{equation} \label{eq.DeltanewDeltaold}
\begin{split}
\Delta_k^{old}(\la) &:=\displaystyle \sum_{i=0}^{k-1} \Delta_{i}(\la), \text{ and }\\ \Delta_k^{new}(\la)& := \Delta_k^{old}(\la)+\Delta_k(\la)=\displaystyle \sum_{i=0}^{k} \Delta_{i}(\la),
\end{split}
\end{equation}
and we will also compute bounds for the Frobenius norm of these error pencils. In our analysis, we will assume that $\delta :=\frac{\|\Delta_S(\la)\|_F}{\|S(\la)\|_F} $ is very small, since in practice is of the order of the machine precision $\epsilon_M$, and we will neglect, when appropriate, terms of order larger than $1$ in $\delta$ to simplify our bounds. Moreover, we will assume that $\delta$ is sufficiently small for guaranteeing that all the steps in the analysis can be performed, for instance, for guaranteeing that some perturbed matrices are invertible, etc. In particular, we have Lemma \ref{cum_errors} for computing bounds of the growth of the cumulative errors $\Delta_k^{new}(\la).$

\begin{lemma}\label{cum_errors}  At each step $k$ of our method, the perturbation $\Delta_k^{new}(\la)$ can be bounded by
	$$  \|\Delta_k^{new}(\la)\|_F \le \sqrt{2}\|\wh S_{k-1}(\la)\|_{2} \|(X_k,Y_k)\|_F +  \|\Delta_k^{old}(\la)\|_F +  {\mathcal O}(\delta^2),
	$$	
	assuming that $\|(X_k,Y_k)\|_F$ is of the order of $\|\Delta_S(\la)\|_F.$
\end{lemma}
\begin{proof}At step $k,$ we have $\wh{S}_{k}(\la)=(I-X_k)\wh{S}_{k-1}(\la)(I-Y_k).$ Therefore
	$$  \Delta^{new}_k(\la) = \Delta^{old}_k(\la) - X_k\widehat S_{k-1}(\la)-\widehat S_{k-1}(\la)Y_k +X_k\widehat S_{k-1}(\la)Y_k.
	$$
	It then follows that the increment (up to ${\mathcal O}(\delta^2)$ terms) is given by
	$$  -X_kS_a-S_aY_k + \la(X_kS_b+S_bY_k) +   {\mathcal O}(\delta^2),$$ where $S_a -\la S_b:=\wh S_{k-1}(\la)$.
	We then use the inequalities
	$$   \|X_kS_a+S_aY_k\|_F^2 \le 2\|S_a\|^2_{2}\|(X_k,Y_k)\|_F ^2, \quad    \|X_kS_b+S_bY_k\|_F^2 \le 2\|S_b\|^2_{2}\|(X_k,Y_k)\|_F ^2
	$$
	and the definition for $\|\wh S_{k-1}(\la)\|_{2}$, to finally get the required bound.
\end{proof}

\subsection{Step 1: Restoring the block anti-triangular structure}\label{restoring_1}
For step 1, that is, restoring the block anti-triangular structure of $S(\la)$ in the perturbed matrix pencil \eqref{eq:perturbed}, we apply a strict equivalence transformation of the type~:
\begin{equation}\label{mat:equiv}
\left[\begin{array}{ccc} I_{(\eta+1)m} & 0 &  0 \\ -X_{21} & I_{\ell} & 0\\ -X_{31}  & -X_{32} & I_{\e n}  \end{array}\right]
\widehat S(\la)
\left[\begin{array}{ccc} I_{(\e+1)n} & -Y_{12} &  -Y_{13} \\ 0 & I_{\ell} & -Y_{23} \\ 0  & 0 & I_{\eta m}  \end{array}\right]
\end{equation}
in order to eliminate the perturbations $\Delta_{23}(\la)$, $\Delta_{32}(\la)$ and $\Delta_{33}(\la)$ of the error matrix pencil $\Delta_0(\la)$. The notation $\widehat S_{ij}^a-\la \widehat S_{ij}^b:=\widehat S_{ij}:= \widehat S_{ij}(\la)$ will be used in this section to refer to sub-blocks of
$\wh S_0(\la)$. Let us write down the equations that we get by setting the blocks (2,3), (3,2) and (3,3) of the matrix in \eqref{mat:equiv} equal to zero~:
\begin{equation}\label{quad}
\begin{split}
\Delta_{23}(\la) :=\Delta_{23}^a-\la \Delta_{23}^b & =  X_{21}\widehat S_{13}+\widehat S_{21}Y_{13}+\widehat S_{22}Y_{23} -X_{21}\widehat S_{11}Y_{13}-X_{21}\widehat S_{12}Y_{23}, \\
\Delta_{32}(\la):=\Delta_{32}^a-\la \Delta_{32}^b & =  \widehat S_{31}Y_{12}+ X_{31}\widehat S_{12}+X_{32} \widehat S_{22} -X_{31}\widehat S_{11}Y_{12}-X_{32}\widehat S_{21}Y_{12}, \\
\Delta_{33}(\la):=\Delta_{33}^a-\la \Delta_{33}^b & =   X_{31}\widehat S_{13}+\widehat S_{31}Y_{13}+X_{32}\Delta_{23}+\Delta_{32}Y_{23} \\
& \quad -X_{31}\widehat S_{11}Y_{13}-X_{32}\widehat S_{21}Y_{13} -X_{31}\widehat S_{12}Y_{23}-X_{32}\widehat S_{22}Y_{23}.
\end{split}
\end{equation}
This is a system of nonlinear matrix equations for the six matrix unknowns $X_{21}, X_{31}, X_{32},$ $Y_{12}, Y_{13}$ and $Y_{23}$. We will show that it is consistent and that it has a solution for which the norms of the unknowns are of the order of $\|\Delta_0 (\la)\|_F$, which implies that there are many terms in the above three equations that are of second order.

Using Kronecker product and the vec$(\cdot)$ notation, the system of matrix equations \eqref{quad} can be rewritten as~:
\begin{equation}\label{eq:quadraticsystem}
\underbrace{\left[\begin{array}{c} \mathrm{vec}(\Delta_{23}^a)\\ \mathrm{vec}(\Delta_{23}^b) \\ \mathrm{vec}(\Delta_{32}^a)\\ \mathrm{vec}(\Delta_{32}^b) \\
	\mathrm{vec}(\Delta_{33}^a)\\ \mathrm{vec}(\Delta_{33}^b) \\ \end{array} \right]}_{:=c} = (T+\Delta T)\underbrace{\left[\begin{array}{c} \mathrm{vec}(X_{21})\\ \mathrm{vec}(Y_{23}) \\ \mathrm{vec}(X_{32})\\ \mathrm{vec}(Y_{12}) \\
	\mathrm{vec}(X_{31})\\ \mathrm{vec}(Y_{13}) \\ \end{array} \right]}_{:=x} - \underbrace{\left[\begin{array}{c} \mathrm{vec}(Z_1)\\ \mathrm{vec}(Z_2) \\ \mathrm{vec}(Z_3)\\ \mathrm{vec}(Z_4) \\
	\mathrm{vec}(Z_5)\\ \mathrm{vec}(Z_6) \\ \end{array} \right]}_{:=z} ,
\end{equation}
where
$$Z_1:=X_{21}\widehat S_{11}^aY_{13}+X_{21}\widehat S_{12}^aY_{23}, \quad Z_2:= X_{21}\widehat S_{11}^bY_{13}+X_{21}\widehat S_{12}^bY_{23},$$       $$Z_3:= X_{31}\widehat S_{11}^aY_{12}+X_{32}\widehat S_{21}^aY_{12}, \quad Z_4:= X_{31}\widehat S_{11}^bY_{12}+X_{32}\widehat S_{21}^bY_{12},$$   $$Z_5:=X_{31}\widehat S_{11}^aY_{13}+X_{32}\widehat S_{21}^aY_{13} +X_{31}\widehat S_{12}^aY_{23}+X_{32}\widehat S_{22}^aY_{23}, $$  $$Z_6:=X_{31}\widehat S_{11}^bY_{13}+X_{32}\widehat S_{21}^bY_{13} +X_{31}\widehat S_{12}^bY_{23}+X_{32}\widehat S_{22}^bY_{23},  $$
\bigskip
{\small $$  \Delta T = \left[\begin{array}{cccccc}
{\Delta_{13}^a}^{T} \otimes I_{\ell} & I_{\eta m}\otimes \Delta_{22}^a  & 0 & 0 & 0 & I_{\eta m}\otimes \Delta_{21}^a \\  {\Delta_{13}^b}^{T} \otimes I_{\ell} & I_{\eta m}\otimes \Delta_{22}^b & 0 & 0 & 0 & I_{\eta m}\otimes \Delta_{21}^b \\
0 & 0  &  {\Delta_{22}^a}^{T} \otimes I_{\e n}  & I_{\ell}\otimes \Delta_{31}^a &  {\Delta_{12}^a}^{T} \otimes I_{\e n} & 0 \\  0 & 0 &  {\Delta_{22}^b}^{T} \otimes I_{\e n}  & I_{\ell}\otimes \Delta_{31}^b &  {\Delta_{12}^b}^{T} \otimes I_{\e n} & 0 \\ 0 & I_{\eta m}\otimes \Delta_{32}^a & {\Delta_{23}^a}^{T} \otimes I_{\e n}  & 0  & {\Delta_{13}^a}^{T} \otimes I_{\e n} & I_{\eta m}\otimes \Delta_{31}^a  \\
0 & I_{\eta m}\otimes \Delta_{32}^b & {\Delta_{23}^b}^{T} \otimes I_{\e n}  & 0  & {\Delta_{13}^b}^{T} \otimes I_{\e n} & I_{\eta m}\otimes \Delta_{31}^b  \end{array} \right],$$
and
$$ T = \left[\begin{array}{cccccc} T_{11} & T_{12} & 0 & 0 & 0 & T_{16} \\  T_{21} & T_{22} & 0 & 0 & 0 & 0 \\
&  &  T_{33} & T_{34} & T_{35} & 0 \\   &  &  T_{43} & T_{44} & 0 & 0 \\ & &  &  & T_{55} & T_{56} \\
& &  &  & T_{65} & T_{66} \end{array} \right],
$$}
with
{\small
	$$\left[\begin{array}{cc} T_{11} & T_{12} \\  T_{21} & T_{22} \end{array} \right]:= \left[ \begin{array}{c|c} E_\eta \otimes I_{m \ell} & I_{\eta m} \otimes A \\ \hline
	F_\eta \otimes I_{m\ell} & I_{\eta m} \otimes I_\ell \end{array} \right], \,
	\left[\begin{array}{cc} T_{33} & T_{34} \\  T_{43} & T_{44} \end{array} \right]:=  \left[ \begin{array}{c|c} A^T\otimes I_{\e n} & I_{\ell} \otimes E_\e\otimes I_n \\ \hline I_\ell \otimes I_{\e n} & I_{\ell} \otimes F_\e\otimes I_n \end{array} \right],$$
	$$\left[\begin{array}{cc} T_{55} & T_{56} \\  T_{65} & T_{66} \end{array} \right]:= \left[ \begin{array}{c|c} E_\eta\otimes I_{m\e n} & I_{\eta m} \otimes E_\e\otimes I_n \\ \hline
	F_\eta \otimes I_{m \e n} & I_{\eta m} \otimes F_\e\otimes I_n \end{array} \right], \quad
	\left\{ \begin{array}{cc} T_{16}  := I_{\eta m} \otimes \mathbf{e}^T_{\e+1}\otimes B \\  T_{35} :=  \mathbf{e}_{\eta+1}^T\otimes C^T \otimes I_{\e n} \end{array} \right. .$$}
We emphasize that the matrices in the two lines above are precisely those appearing in equations \eqref{singval_2}, \eqref{singval_1} and \eqref{singval_3}, respectively.

The smallest singular value of $T$ and the 2--norm of $\Delta T$ will be needed in the analysis of the bound for the structured backward errors. More precisely for proving that \eqref{quad} is consistent and bounding the norm of one of its solutions.
A  lower bound for $\sigma_{\min} (T)$ and an upper bound for $\|\Delta T\|_2$ are given in Lemma \ref{lem:sigmamin} and Lemma \ref{lem:deltat}, respectively.

\begin{lemma}\label{lem:sigmamin}
	Let $T$ be the matrix in \eqref{eq:quadraticsystem}. Let $\alpha:=1+2\e\max(1,\|A\|_2^\e),$ $\beta:=1+2\eta\max(1,\|A\|_2^\eta),$ $\gamma:=\frac{\e+\eta}{2\sqrt{2}}$ and $s:= \max(\alpha,\beta,\gamma)+\gamma(\beta\|B\|_2+\alpha\|C\|_2)$ then
	$$ \sigma_{\min}(T) \ge \frac{1}{s}.
	$$
\end{lemma}
\begin{proof}
	If we partition the matrix $T$ as a block triangular matrix
	$$ T=\left[\begin{array}{ccc} T_{1} & 0 & T_B \\   & T_{2} & T_C  \\  &  &  T_{3}  \end{array} \right] ,  $$
	then the diagonal blocks have full row ranks because their smallest singular values are strictly larger than zero according to Lemmas \ref{lem:singval_2}, \ref{lem:singval_1} and \ref{lem:singval_3}, respectively. Therefore, they are
	right invertible, with Moore--Penrose pseudoinverses $T_i^r$ satisfying $T_i T_i^r = I$, for $i=1,2,3$. Moreover, $\|T_1^r\|_2 = \omega_2^{-1}$,  $\|T_2^r\|_2 = \omega_1^{-1}$ and $\|T_3^r\|_2 = \omega_3^{-1}$, with $\omega_1, \omega_2$ and $\omega_3$ as in \eqref{singval_1}, \eqref{singval_2} and \eqref{singval_3}. A right inverse $T^r$ for $T$ is given by
	$$T^r=\left[\begin{array}{ccc} T_{1}^r & 0 & -T_{1}^rT_BT_{3}^r \\   & T_{2}^r & -T_{2}^rT_CT_{3}^r  \\  &  &  T_{3}^r  \end{array} \right] $$
	since $T T^r=I$. It then follows that the smallest singular value of $T$ is lower bounded by $\|T^r\|_2^{-1}$. This right inverse can be written as the sum of three matrices (one of them being  $\mathrm{diag} (T_1^r, T_2^r, T_3^r)$), and the 2-norm of each of them can be upper bounded  using the results of Section \ref{SVbounds} and the fact that $\|T_B\|_2=\|B\|_2$ and  $\|T_C\|_2=\|C\|_2$. We then obtain the bound~:
	\begin{equation*}
	\begin{split}
	\sigma_{\min}(T)  \ge & 1/\left[ \max(\omega_1^{-1},\omega_2^{-1},\omega_3^{-1})+\omega_3^{-1}(\omega_2^{-1}\|B\|_2+\omega_1^{-1}\|C\|_2)\right] \\
	\ge & 1/\left[\max(\alpha,\beta,\gamma)+\gamma(\beta\|B\|_2+\alpha\|C\|_2)\right],
	\end{split}
	\end{equation*}
	by taking into account inequalities \eqref{eq:omega1}, \eqref{eq:omega2}, \eqref{eq:omega3}.
\end{proof}

\begin{lemma}\label{lem:deltat}
	Let $\Delta T$ be the matrix in \eqref{eq:quadraticsystem} and let $\Delta_{S}(\la)$ be the pencil in \eqref{eq:perturbed}. Then $$\| \Delta T \|_{2} \leq  \sqrt{3} \| \Delta_{S}(\la) \|_{2}. $$
\end{lemma}
\begin{proof} We consider a permutation matrix $P$ such that
{\small	\begin{equation*}
	\begin{split}
	\Delta T  & =   \left[\begin{array}{cc|cc|cc}
	{\Delta_{13}^a}^{T} \otimes I_{\ell} & 0 & 0 & 0 & I_{\eta m}\otimes \Delta_{22}^a   & I_{\eta m}\otimes \Delta_{21}^a \\  {\Delta_{13}^b}^{T} \otimes I_{\ell} & 0 & 0 & 0 & I_{\eta m}\otimes \Delta_{22}^b & I_{\eta m}\otimes \Delta_{21}^b \\
	0 & I_{\ell}\otimes \Delta_{31}^a &  {\Delta_{22}^a}^{T} \otimes I_{\e n}  & {\Delta_{12}^a}^{T} \otimes I_{\e n}  & 0  & 0 \\  0 &
	I_{\ell}\otimes \Delta_{31}^b &  {\Delta_{22}^b}^{T} \otimes I_{\e n}  & {\Delta_{12}^b}^{T} \otimes I_{\e n}  &  0 & 0 \\ 0 & 0 & {\Delta_{23}^a}^{T} \otimes I_{\e n}  & {\Delta_{13}^a}^{T} \otimes I_{\e n}  &  I_{\eta m}\otimes \Delta_{32}^a & I_{\eta m}\otimes \Delta_{31}^a  \\
	0 & 0 & {\Delta_{23}^b}^{T} \otimes I_{\e n}  & {\Delta_{13}^b}^{T} \otimes I_{\e n}  &  I_{\eta m}\otimes \Delta_{32}^b & I_{\eta m}\otimes \Delta_{31}^b  \end{array} \right]P \\   & :=   \left[T_1| T_2|T_3\right]P.
	\end{split}
	\end{equation*}}
	Using properties of norms and Kronecker products (see \cite[Chapter 4]{HyJ}) we have that $\| T_i \|_{2} \leq  \| \Delta_{S}(\la) \|_{2}$ for $i=1,2,3.$ Finally, by \cite[Lemma 3.5]{HLT},
	$$\| \Delta T \|_{2} \leq \sqrt{3}\max \{\| T_1 \|_{2},\| T_2 \|_{2},\| T_3 \|_{2}\} \leq   \sqrt{3} \| \Delta_{S}(\la) \|_{2}. $$	 \end{proof}

In order to prove that the system of nonlinear matrix equations \eqref{quad} is consistent, first, we remove quadratic terms in $X_{ij}$ and $Y_{ij}$ of these equations and we get the following system of linear equations~:
\begin{equation*}
\begin{split}
\Delta_{23}(\la) & =  X_{21}\widehat S_{13}+\widehat S_{21}Y_{13}+\widehat S_{22}Y_{23} , \\
\Delta_{32}(\la) & =  \widehat S_{31}Y_{12}+ X_{31}\widehat S_{12}+X_{32} \widehat S_{22}, \\
\Delta_{33}(\la) & =   X_{31}\widehat S_{13}+\widehat S_{31}Y_{13}+X_{32}\Delta_{23}+\Delta_{32}Y_{23}.
\end{split}
\end{equation*}
This linear system of matrix equations can be rewritten as the underdetermined linear system~:
\begin{equation}\label{eq:linearsystem}
(T+\Delta T) x  =  c ,
\end{equation}
with the same notation as in \eqref{eq:quadraticsystem}. Next we prove that \eqref{eq:linearsystem} is consistent for any right hand side if $\Delta T$ is sufficiently small. From the minimum norm solution of \eqref{eq:linearsystem}, we obtain in Theorem \ref{theo:bound} that there exists a solution for the quadratic system \eqref{eq:quadraticsystem} under certain conditions and bound its norm.

\begin{lemma}\label{lem:minsolution} Let $(T+\Delta T)x=c$ be the underdetermined linear system in \eqref{eq:linearsystem}, and let us assume that $\sigma_{\min}(T)>\|\Delta T\|_2.$ Then $(T+\Delta T)x=c$ is consistent and its minimum norm solution $(X^0,Y^0):=(X_{21}^0,X_{31}^0,X_{32}^0,Y_{12}^0,Y_{13}^0,Y_{23}^0)$ satisfies
	$$\|(X^0,Y^0)\|_F \leq \frac{1}{\sigma}\|(\Delta_{23}(\la),\Delta_{32}(\la),\Delta_{33}(\la))\|_F,$$
	where $\sigma:= \sigma_{\min}(T)-\|\Delta T\|_2.$
\end{lemma}
\begin{proof}Analogous proof as for \cite[Lemma 5.6]{DLPV}.
\end{proof}

The notation $\sigma:= \sigma_{\min}(T)-\|\Delta T\|_2$ has been chosen to remind that $\sigma$ is a lower bound for the smallest singular value of $T+\Delta T,$ since $\sigma_{\min}(T+\Delta T) \geq \sigma_{\min}(T)-\|\Delta T\|_2$ by Weyl's perturbation theorem for singular values \cite[Theorem 3.3.16]{HyJ}. Lemma \ref{lem:delta} gives a sufficient condition on $\| \Delta_{S}(\la) \|_2$ that guarantees $\sigma >0$ and, hence, that allows us to apply Lemma \ref{lem:minsolution}.

\begin{lemma}\label{lem:delta} Consider the real number $s$ defined as in Lemma \ref{lem:sigmamin}. Let $T$ and $\Delta T$ be the matrices in \eqref{eq:linearsystem}, and let $\Delta_{S}(\la)$ be the pencil in \eqref{eq:perturbed}. If $\|\Delta_{S}(\la) \|_2 < \frac{1}{2s} $ then
	$$\sigma = \sigma_{\min}(T)-\| \Delta T \|_{2} >  \frac{2-\sqrt{3}}{2s}  >0.$$
\end{lemma}
\begin{proof} If $\| \Delta_{S}(\la) \|_2 < \frac{1}{2s} $ we have, by Lemmas \ref{lem:sigmamin} and \ref{lem:deltat}, that $\sigma_{\min}(T)-\| \Delta T \|_{2} \geq \frac{1}{s} - \sqrt{3}\, \|\Delta_{S}(\la) \|_2 > \frac{2-\sqrt{3}}{2s}  >0.$
\end{proof}

Theorem \ref{theo:bound} establishes conditions in order the system of matrix equations \eqref{quad} to have a solution as we announced. Moreover, it gives an upper bound for the Frobenius norm of this solution. We remark that Theorem \ref{theo:bound} is similar to \cite[Theorem 5.8]{DLPV}, though the involved systems of matrix equations are very different from each other. Therefore, some details in the proof of Theorem \ref{theo:bound} are omitted since can be found in \cite{DLPV}.

\begin{theorem}\label{theo:bound} There exists a solution $(X,Y):=(X_{21},X_{31},X_{32},Y_{12},Y_{13},Y_{23})$ of the quadratic system of equations \eqref{eq:quadraticsystem} satisfying
	$$\|(X,Y)\|_F \leq 2\frac{\theta}{\sigma},$$
	whenever
	\begin{equation}\label{condition}
	\sigma >0 \quad\mbox{ and }\quad \frac{\theta \omega}{\sigma^2} < \frac{1}{4},
	\end{equation}
	where $\omega := \| (M(\la),A-\la I_{\ell},B, C)\|_F + \| \Delta_{S}(\la) \|_F,$ $\theta := \|(\Delta_{23}(\la),\Delta_{32}(\la),\Delta_{33}(\la))\|_F,$ and $\sigma =\sigma_{\min}(T)-\|\Delta T\|_2.$
\end{theorem}
\begin{proof} Since $\sigma>0,$ we can apply Lemma \ref{lem:minsolution} and consider $(X^0,Y^0)$ the minimum norm solution of \eqref{eq:linearsystem}. Let $$x_{0}:=\left[\begin{array}{cccccc} \mathrm{vec}(X_{21}^0)^{T} \, \mathrm{vec}(Y_{23}^0)^{T} \, \mathrm{vec}(X_{32}^0)^{T} \, \mathrm{vec}(Y_{12}^0)^{T} \,
	\mathrm{vec}(X_{31}^0)^{T} \, \mathrm{vec}(Y_{13}^0)^{T}  \end{array} \right]^{T}.$$ Let us define the sequence $\{(X^i,Y^i):=(X_{21}^i,X_{31}^i,X_{32}^i,Y_{12}^i,Y_{13}^i,Y_{23}^i)\}_{i=0}^{\infty}$ such that, for each $i>0,$ $(X^i,Y^i)$ is the minimum norm solution of the linear system
	\begin{equation}\label{eq:system2}
	(T+\Delta T) \left[\begin{array}{c} \mathrm{vec}(X_{21}^i)\\ \mathrm{vec}(Y_{23}^i) \\ \mathrm{vec}(X_{32}^i)\\ \mathrm{vec}(Y_{12}^i) \\
	\mathrm{vec}(X_{31}^i)\\ \mathrm{vec}(Y_{13}^i) \\ \end{array} \right]= c +
	\left[\begin{array}{c} \mathrm{vec}(Z_1^{i-1})\\ \mathrm{vec}(Z_2^{i-1}) \\ \mathrm{vec}(Z_3^{i-1})\\ \mathrm{vec}(Z_4^{i-1}) \\
	\mathrm{vec}(Z_5^{i-1})\\ \mathrm{vec}(Z_6^{i-1}) \\ \end{array} \right] ,
	\end{equation}
	where
	$$Z_1^{i-1}:=X_{21}^{i-1}\widehat S_{11}^aY_{13}^{i-1}+X_{21}^{i-1}\widehat S_{12}^aY_{23}^{i-1},\quad Z_2^{i-1}:= X_{21}^{i-1}\widehat S_{11}^bY_{13}^{i-1}+X_{21}^{i-1}\widehat S_{12}^bY_{23}^{i-1},$$       $$Z_3^{i-1}:= X_{31}^{i-1}\widehat S_{11}^aY_{12}^{i-1}+X_{32}^{i-1}\widehat S_{21}^aY_{12}^{i-1},\quad Z_4^{i-1}:= X_{31}^{i-1}\widehat S_{11}^bY_{12}^{i-1}+X_{32}^{i-1}\widehat S_{21}^bY_{12}^{i-1},$$   $$Z_5^{i-1}:=X_{31}^{i-1}\widehat S_{11}^aY_{13}^{i-1}+X_{32}^{i-1}\widehat S_{21}^aY_{13}^{i-1}+X_{31}^{i-1}\widehat S_{12}^aY_{23}^{i-1}+X_{32}^{i-1}\widehat S_{22}^aY_{23}^{i-1}, \quad \text{and}$$  $$Z_6^{i-1}:=X_{31}^{i-1}\widehat S_{11}^bY_{13}^{i-1}+X_{32}^{i-1}\widehat S_{21}^bY_{13}^{i-1} +X_{31}^{i-1}\widehat S_{12}^bY_{23}^{i-1}+X_{32}^{i-1}\widehat S_{22}^bY_{23}^{i-1}. $$
	Note that the minimum norm solution of \eqref{eq:system2} is obtained by multiplying the right hand side of \eqref{eq:system2} by the Moore-Penrose pseudoinverse of $T + \Delta T$, denoted by $(T + \Delta T)^\dagger$, and that $x_0 = (T + \Delta T)^\dagger c$.

	Now we assume that $\frac{\theta \omega}{\sigma^2} < \frac{1}{4}$ holds. Then we can prove that the sequence $\{(X^i,Y^i)\}_{i=0}^{\infty}$ converges to a solution $(X,Y)$ of the quadratic system of equations \eqref{eq:quadraticsystem} analogously as it is done in \cite[Theorem 5.8]{DLPV}. For that, we have to take into account that, if $\|(X^{i-1},Y^{i-1})\|_{F}\leq \rho_{i-1},$ then {\small
		\begin{equation*}
		\begin{split}
	&	\|(X^i, Y^i)\|_F  \\ &  \leq  \|(X^0,Y^0)\|_F + \|(T+\Delta T)^{\dagger}\|_2 \left\| \left[\begin{array}{cc}
		X_{21}^{i-1} & 0 \\
		X_{31}^{i-1} & X_{32}^{i-1}
		\end{array}\right]\left[\begin{array}{cc}
		\widehat S_{11} & \widehat S_{12} \\
		\widehat S_{21} & \widehat S_{22}
		\end{array}\right]\left[\begin{array}{cc}
		Y_{12}^{i-1} & Y_{13}^{i-1} \\
		0 & Y_{23}^{i-1}
		\end{array}\right] \right\|_F \\
		 & \leq \rho_{0}+\sigma^{-1}\rho_{i-1}^{2}\omega :=\rho_i,
		\end{split}
		\end{equation*}}
	where $\|(X^0,Y^0)\|_F \leq \theta\sigma^{-1}:=\rho_0.$ Therefore, we can define the same fixed point iteration as in the proof of \cite[Theorem 5.8]{DLPV} and we obtain that the sequence is bounded, i.e., $\|(X^i,Y^i)\|_{F}\leq \rho,$ with $\rho<2\sigma^{-1}\theta$, for all $i \geq 0$.  In addition, if we define the sequence $\{C_{i}:=(X^{i+1},Y^{i+1})-(X^i,Y^i)\}_{i=0}^{\infty}$ then
{\small	\begin{equation*}
	\begin{split}
	\|C_i\|_F \leq & \|(T+\Delta T)^{\dagger}\|_2 \left(\left\| \left[\begin{array}{cc}
	X_{21}^{i} & 0 \\
	X_{31}^{i}& X_{32}^{i}
	\end{array}\right]\left[\begin{array}{cc}
	\widehat S_{11} & \widehat S_{12} \\
	\widehat S_{21} & \widehat S_{22}
	\end{array}\right]\left[\begin{array}{cc}
	Y_{12}^{i} & Y_{13}^{i} \\
	0 & Y_{23}^{i}
	\end{array}\right] \right. \right.  \\ & \left. \left. -  \left[\begin{array}{cc}
	X_{21}^{i-1} & 0 \\
	X_{31}^{i-1} & X_{32}^{i-1}
	\end{array}\right]\left[\begin{array}{cc}
	\widehat S_{11} & \widehat S_{12} \\
	\widehat S_{21} & \widehat S_{22}
	\end{array}\right]\left[\begin{array}{cc}
	Y_{12}^{i-1}& Y_{13}^{i-1} \\
	0 & Y_{23}^{i-1}
	\end{array}\right]  \right\|_F \right) \\
	\leq & \|(T+\Delta T)^{\dagger}\|_2 \left(\left\| \left[\begin{array}{cc}
	X_{21}^{i}-X_{21}^{i-1} & 0 \\
	X_{31}^{i}-X_{31}^{i-1}& X_{32}^{i}-X_{32}^{i-1}
	\end{array}\right]\left[\begin{array}{cc}
	\widehat S_{11} & \widehat S_{12} \\
	\widehat S_{21} & \widehat S_{22}
	\end{array}\right]\left[\begin{array}{cc}
	Y_{12}^{i} & Y_{13}^{i} \\
	0 & Y_{23}^{i}
	\end{array}\right] \right\|_F  \right. \\ &  \left. +   \left\| \left[\begin{array}{cc}
	X_{21}^{i-1} & 0 \\
	X_{31}^{i-1} & X_{32}^{i-1}
	\end{array}\right]\left[\begin{array}{cc}
	\widehat S_{11} & \widehat S_{12} \\
	\widehat S_{21} & \widehat S_{22}
	\end{array}\right]\left[\begin{array}{cc}
	Y_{12}^{i}-Y_{12}^{i-1}& Y_{13}^{i}-Y_{13}^{i-1} \\
	0 & Y_{23}^{i}-Y_{23}^{i-1}
	\end{array}\right]  \right\|_F \right) \\
	\leq & 2 \sigma^{-1}\rho\omega\|C_{i-1}\|_F.
	\end{split}
	\end{equation*}}
	The above inequality implies that $\{(X^i,Y^i)\}_{i=0}^{\infty}$ is a Cauchy sequence, since $ 2 \sigma^{-1}\rho\omega<1.$ Thus, taking limits in both sides of \eqref{eq:system2}, we see that $\{(X^i,Y^i)\}_{i=0}^{\infty}$ converges to a solution $(X,Y)$ of the system of equations in \eqref{eq:quadraticsystem} with $\|(X,Y)\|_{F}\leq \rho.$
\end{proof}

Theorem \ref{theo:bound}, together with Lemma \ref{lem:delta}, allow us to prove in Theorem \ref{theo:bound2} that there exists a solution $(X,Y)$ of \eqref{quad} which is of the order of the perturbation $\Delta_{S}(\la)$ whenever $\|\Delta_{S}(\la)\|_F$ is properly upper bounded.

\begin{theorem}\label{theo:bound2} Consider the real number $s$ defined as in Lemma \ref{lem:sigmamin}. Let $S(\la)$ be a block Kronecker linearization as in \eqref{eq:pencil}, and let $\Delta_{S}(\la)$ be a perturbation of $S(\la)$ as in \eqref{eq:perturbed} such that
	\begin{equation}\label{norm_pert}
	\|\Delta_S(\la)\|_F < \left( \frac{2-\sqrt{3}}{4s}\right)^2\frac{1}{1 + \| (M(\la),A-\la I_{\ell},B, C )\|_F }.
	\end{equation}
	Then there exists a solution $(X,Y):= (X_{21},X_{31},X_{32},Y_{12},Y_{13},Y_{23}) $ of the quadratic system of matrix equations in \eqref{quad} that satisfies
	\begin{equation}
	\|(X,Y)\|_F \leq \frac{4s\|\Delta_S(\la)\|_F}{2-\sqrt{3}}.
	\end{equation}
\end{theorem}
\begin{proof} We have	\begin{equation*}
	\|\Delta_S(\la)\|_F < \left( \frac{2-\sqrt{3}}{4s}\right)^2\frac{1}{1 + \| (M(\la),A-\la I_{\ell},B, C )\|_F }\leq \frac{1}{2s}
	\end{equation*}
	since $s\geq 1.$ Then, by Lemma \ref{lem:delta}, $\sigma = \sigma_{\min}(T)-\|\Delta T\|_2  >  \frac{2-\sqrt{3}}{2s}  >0.$ In addition, using the same notation as in Theorem \ref{theo:bound},
	\begin{equation*}
	\frac{\theta \omega}{\sigma^2} \leq \dfrac{\|\Delta_S (\la)\|_F (\| (M(\la),A-\la I_{\ell},B,C)\|_F + \|\Delta_{S}(\la) \|_F)}{\left(\frac{2-\sqrt{3}}{2s}\right)^2} <  \frac{1}{4},
	\end{equation*}
	by \eqref{norm_pert}. Therefore, conditions in \eqref{condition} hold and, by Theorem \ref{theo:bound}, there exists a solution $(X,Y)$ of the system in \eqref{quad} satisfying
	$$\|(X,Y)\|_{F} \leq 2\frac{\theta}{\sigma}\leq \frac{4s \|\Delta_S (\la)\|_F }{2-\sqrt{3}}.$$
\end{proof}

After restoring the block anti-triangular structure of $S(\la)$, we get the perturbation error $\Delta^{new}_1(\la)$ defined in \eqref{eq.DeltanewDeltaold}. The following first order bound for the norm of $\Delta^{new}_1(\la)$ in Corollary \ref{C1} follows from Lemma \ref{cum_errors} and Theorem \ref{theo:bound2}.

\begin{corollary}  \label{C1}   Let us define the scalar $f_1:=\frac{4\sqrt{2}s}{2-\sqrt{3}}$. Then
\begin{equation*}
\begin{split}
	\|\Delta^{new}_1(\la)\|_F & \le [1 + f_1 \|  \widehat S_0 (\la)\|_2] \, \|\Delta_S(\la)\|_F +{\mathcal O}(\delta^2) \\ 
 & \le [1 + f_1 \|  S(\la)\|_2] \, \|\Delta_S(\la)\|_F +{\mathcal O}(\delta^2) .
\end{split}
\end{equation*}

\end{corollary}

\subsection{Step 2: Restoring the Kronecker blocks $K_1(\la)$, $K_2(\la)$ and the identity $I_\ell$}\label{restoring_2}

At this stage we have obtained a pencil $ \widehat S_1(\la)=  S(\la)+ \Delta^{new}_1(\la)$  of the type
\begin{equation} \label{eqstep2}
\widehat S_1(\la):=\left[\begin{array}{ccc} \widehat M(\la) & \widehat C(\la) &  \widehat K_2^T(\la)  \\ \widehat B(\la) & \widehat A- \la \widehat I_\ell & 0\\
\widehat K_1(\la) &   0 &  0
\end{array}\right],
\end{equation}
where the zero blocks below the anti-diagonal are exact and $\widehat{S}_1 (\la)$ is strictly equivalent to $\widehat{S} (\la)$.
In this subsection, we will use $\Delta^{a}_{ij}-\la \Delta^{b}_{ij}$ to denote the corresponding blocks of the updated perturbation matrix
$\Delta^{new}_1(\la)$. We assume that the norm of the perturbation $\Delta^{new}_1(\la)$ is small enough for $ \widehat K_1(\la) $ and $ \widehat K_2(\la) $ to be also minimal bases with row degrees all equal to $1$ and the row degrees of their dual minimal bases all equal to $\epsilon$ and $\eta$, respectively \cite[Corollary 5.15]{DLPV}. Thus, $ \widehat K_1(\la) $ and $ \widehat K_2(\la) $ have the same Kronecker canonical forms as $K_1(\la)$ and $K_2(\la)$, respectively, and are strictly equivalent to them. We will then perform step 2, that is, an updating block-diagonal strict equivalent transformation of the type
{\small \begin{equation} \label{trstep2}
	\left[\begin{array}{ccc} I_{(\eta+1)m}-X_{11} & 0 &  0 \\  0 & I_{\ell}-X_{22} & 0\\ 0 & 0 & I_{\e n}-X_{33} \end{array}\right]
	\widehat S_1(\la)  \left[\begin{array}{ccc} I_{(\epsilon + 1)n}-Y_{11} & 0 & 0 \\ 0 & I_{\ell}-Y_{22} & 0 \\ 0  & 0 & I_{\eta m}-Y_{33}  \end{array}\right]
	\end{equation}}
such that
$$   (I-X_{33})\widehat K_1(\la) (I-Y_{11}) = K_1(\la), \quad (I-X_{11})\widehat K_2^T(\la) (I-Y_{33}) = K_2^T(\la),
$$
and $$ (I-X_{22})\widehat I _\ell (I-Y_{22}) = I_\ell. $$ In the last three equations the sizes of some identity matrices are not specified for simplicity.
Clearly, these three problems are independent from each other and can be treated separately.

Let us first look at the equation restoring $K_1(\la)$.
As pointed out in Section \ref{Sylvester}, this can be reduced to the solution of a Sylvester equation. Let
$$ \widehat K_1(\la) = K_1(\la) + \Delta_{K_1}(\la) := L_\e(\la)\otimes I_n +  \Delta_{K_1}(\la):= (E_\e- \la F_\e)\otimes I_n + (\Delta^a_{31}-\la \Delta^b_{31}).$$
Then, making the change of variables $Y_{11}:=Y$ and $X_{33}:=X(I+X)^{-1}$, it suffices to solve
$$  (K_1(\la)+\Delta_{K_1}(\la))Y+ X K_1(\la) = \Delta_{K_1}(\la),
$$
or, equivalently,
\begin{equation} \label{eq.auxfro1}
\left[ \begin{array}{c|c} E_\e^T \otimes I_{n\e n} & I_{(\e+1)n} \otimes ( E_\e\otimes I_n + \Delta^a_{31}) \\ \hline
F_\e^T \otimes I_{n\e n} & I_{(\e+1)n} \otimes (F_\e\otimes I_n + \Delta^b_{31}) \end{array} \right]
\left[ \begin{array}{c} \mathrm{vec}(X) \\ \mathrm{vec}(Y) \end{array} \right] =
\left[ \begin{array}{c} \mathrm{vec}(\Delta^a_{31}) \\ \mathrm{vec}(\Delta^b_{31}) \end{array} \right] .
\end{equation}
By Lemma \ref{lem:singval_4}, the smallest singular value of the unperturbed problem
satisfies
$$ \sigma_{2 \e n(\e+1)n}\left[ \begin{array}{c|c} E_\e^T \otimes I_{n\e n} & I_{(\e+1)n} \otimes E_\e\otimes I_n \\ \hline
F_\e^T \otimes I_{n\e n} & I_{(\e+1)n} \otimes F_\e\otimes I_n \end{array} \right] \ge \frac{3}{4\e-1}. $$
Then, by using Weyl's perturbation theorem for singular values \cite[Theorem 3.3.16]{HyJ}, one obtains the following bound for the minimum norm solution of \eqref{eq.auxfro1}
$$  \|(X,Y)\|_F  \le  \left[\frac{3}{4\e-1} - \|\Delta^a_{31}\|_2 - \|\Delta^b_{31}\|_2\right]^{-1}\|(\Delta^a_{31},\Delta^b_{31})\|_F , $$
assuming that the perturbation is small enough for satisfying $\frac{3}{4\e-1} - \|\Delta^a_{31}\|_2 - \|\Delta^b_{31}\|_2 > 0$.
In addition,
$$  \|(X_{33},Y_{11})\|_F \le  \|(X,Y)\|_F/(1 - \|(X,Y)\|_F).
$$
Since $\| \Delta^a_{31}\|_2$ and $\| \Delta^b_{31}\|_2$ are of the order of $\delta$, finally yields
\begin{equation} \label{ineq1}    \|(X_{33},Y_{11})\|_F  \le   \frac{4\e-1}{3}\|(\Delta^a_{31},\Delta^b_{31})\|_F + {\mathcal O}(\delta^2),
\end{equation} \label{th:2}
by neglecting quantities of the order of ${\mathcal O}(\delta^2)$.

The problem for restoring $K_2(\la)$ is clearly dual to the problem of $K_1(\la)$ and will therefore yield the bound
\begin{equation} \label{ineq2}    \|(X_{11},Y_{33})\|_F  \le   \frac{4\eta-1}{3}\|(\Delta^a_{13},\Delta^b_{13})\|_F + {\mathcal O}(\delta^2).
\end{equation}

The problem of restoring $I_\ell$ amounts to solving $(I_\ell -X_{22})(I_\ell +\Delta^b_{22})(I_\ell -Y_{22})=I_\ell$, with $\widehat{I}_\ell = I_\ell + \Delta_{22}^b$. There are many possible solutions. A very simple one is to take $Y_{22} = 0$ and $I_\ell -X_{22} = (I_\ell +\Delta^b_{22})^{-1}$, assuming $\Delta^b_{22}$ is small enough for the inverse to exist. This means that $X_{22} = \Delta^b_{22} + {\mathcal O} (\|\Delta^b_{22} \|_F^2)$ and
\begin{equation} \label{ineq3}  \|(X_{22},Y_{22})\|_F  =  \|\Delta^b_{22}\|_F + {\mathcal O}(\delta^2).
\end{equation}


We summarize this discussion in the following Theorem.
\begin{theorem}\label{theo:bound3}
	Let the pencil $\widehat S_1(\la)$ have the block anti-triangular form given in \eqref{eqstep2}. If $\max (\e ,\eta) > 0$, then the updating strict equivalence transformation $(I-X)\widehat S_1(\la)(I-Y)$
	detailed in \eqref{trstep2} exists and can be bounded by
	$$ \|(X,Y)\|_F \le \frac{4\max(\e,\eta)-1}{3} \| \Delta^{new}_1(\la)\|_F + {\mathcal O}(\delta^2). $$
\end{theorem}
\begin{proof}
	The bound for $\|(X,Y)\|_F$ follows directly from the identity
	$$\|(X,Y)\|_F^2=\|(X_{11},Y_{33})\|_F^2+\|(X_{22},Y_{22})\|_F^2+\|(X_{33},Y_{11})\|_F^2,$$
	from the inequality
	$$\|(\Delta^a_{13},\Delta^b_{13})\|_F^2+\|\Delta^b_{22}\|_F^2+\|(\Delta^a_{31},\Delta^b_{31})\|_F^2 \le \| \Delta^{new}_1(\la)\|_F^2 $$ and from the individual inequalities \eqref{ineq1}, \eqref{ineq2} and \eqref{ineq3}.
\end{proof}

The following first order bound in Corollary \ref{C2} for the norm of the perturbation error $\Delta^{new}_2(\la)$ follows from Lemma \ref{cum_errors}, Theorem \ref{theo:bound3} and Corollary \ref{C1}.

\begin{corollary}   \label{C2}  Let us define the scalar $f_2:=\frac{\sqrt{2}(4\max(\e,\eta)-1)}{3}$. Then
\begin{equation*}
\begin{split}
	\|\Delta^{new}_2(\la)\|_F &  \le [1 + f_2 \| \widehat S_1(\la)\|_2] \, \|\Delta_1^{new}(\la)\|_F +{\mathcal O}(\delta^2) \\
	& \le [1 + f_2 \|  S(\la)\|_2]\, \|\Delta_1^{new}(\la)\|_F +{\mathcal O}(\delta^2) .
	\end{split}
\end{equation*}

\end{corollary}

\subsection{Step 3: Restoring the constant $B$ and $C$ matrices}\label{restoring_3}
From steps 1 and 2, described in the previous subsections, we have obtained a pencil $ \widehat S_2(\la)=  S(\la)+ \Delta^{new}_2(\la)$  of the type
\begin{equation} \label{eqstep3} \widehat S_2(\la):=\left[\begin{array}{ccc} \widehat M(\la) & \widehat C(\la) &  K_2^T(\la)  \\ \widehat B(\la) & \widehat A- \la I_\ell & 0\\
K_1(\la) &   0 &  0
\end{array}\right]
\end{equation}
strictly equivalent to $\widehat{S} (\la)$. We emphasize that the blocks $\widehat{M} (\la), \widehat{B} (\la), \widehat{C} (\la)$ and the matrix $\widehat{A}$ are obviously different in \eqref{eqstep3} and in \eqref{eqstep2}. We use the same symbols for avoiding a cumbersome notation.
In this subsection, we will use $\Delta_{ij} (\la) =  \Delta^{a}_{ij}-\la \Delta^{b}_{ij}$ to denote the corresponding blocks of the updated perturbation matrix
$\Delta^{new}_2(\la)$.
In this third step, we will restore the pencil $\widehat S_2(\la)$ to one where the blocks
$$ \widehat B(\la)=B \widehat K_1 + \Delta_{21}(\la) , \quad  \mathrm{and} \quad \widehat C(\la)= \widehat K_2^T C + \Delta_{12}(\la) $$
are transformed to $\widetilde B \widehat K_1$ and $\widehat K_2^T \widetilde C$, respectively. We recall that
$$
K_1(\la)= L_\epsilon(\la)\otimes I_n,\quad \widehat K_1=\mathbf{e}_{\epsilon+1}^T\otimes I_n, \quad K_2(\la)= L_\eta(\la)\otimes I_m, \quad \widehat K_2=\mathbf{e}_{\eta+1}^T\otimes I_m,
$$
where $\mathbf{e}_k$ is the standard $k$th unit vector of dimension $k$ and $L_k(\la)$ is the classical Kronecker block of dimension $k\times(k+1)$, as introduced below \eqref{block_kronecker_lin}.
We will construct for this a strict equivalence transformation of the type
{\small \begin{equation}\label{trstep3}
	\begin{split} & \left[\begin{array}{ccc} I_{m(\eta+1)} & -X_{12}  & 0  \\  & I_\ell & -X_{23} \\   &  &  I_{n\e} \end{array}\right] \widehat S_2(\la)\left[\begin{array}{ccc}  I_{n(\e+1)} & & \\ -Y_{21} & I_\ell & \\ 0  & -Y_{32} &  I_{m\eta} \end{array}\right] \\  = & \left[\begin{array}{ccc} \widetilde M(\la) &  \widehat K_2^T\widetilde C &  K_2^T(\la)  \\ \widetilde B \widehat K_1 & \widehat A- \la I_\ell & 0\\
	K_1(\la) &   0 &  0
	\end{array}\right] \end{split} \end{equation}}
The problems for $\widehat B(\la)$ and $\widehat C(\la)$ can again be treated separately.
Let us first focus on the subsystem
{\small \begin{equation*}
	\begin{split} & \left[\begin{array}{ccc} I_{m(\eta+1)} & -X_{12} \\ & I_{\ell} \end{array} \right] \left[\begin{array}{ccc} \widehat C(\la) & L_\eta^T(\la)\otimes I_m  \\ \widehat A-\la I_\ell  & 0
	\end{array}\right] \left[\begin{array}{ccc} I_{\ell} & \\  -Y_{32} & I_{m\eta} \end{array} \right] \\  =  & \left[\begin{array}{ccc}   \mathbf{e}_{\eta+1} \otimes \widetilde C  &   L_\eta^T(\la)\otimes I_m  \\ \widehat A -\la I_\ell & 0
	\end{array}\right].
	 \end{split} \end{equation*}}
If we partition the matrices $X_{12}$, $Y_{32}$ and $\widehat C(\la)$ as follows :
$$  X_{12} :=\left[\begin{array}{c} E_1 \\ E_2 \\ \vdots \\ E_\eta \\ E_{\eta+1 } \end{array} \right] , \quad
Y_{32}:= \left[\begin{array}{c} F_1 \\ F_2 \\ \vdots \\ F_\eta  \end{array}\right], \quad   \widehat C(\la)  := \left[\begin{array}{c} C_{01} \\ C_{02} \\ \vdots \\ C_{0\eta} \\ C_{0(\eta+1)} \end{array} \right]  - \left[\begin{array}{c} C_{11} \\ C_{12} \\ \vdots \\ C_{1\eta} \\ C_{1(\eta+1)} \end{array} \right] \la,
$$
where all blocks have dimension $m\times \ell$, then we need to solve the following system of equations
\begin{equation*}
\begin{split}
 & \left[\begin{array}{cccccccc} E_1 & F_1 & E_2  & ... & F_\eta & E_{\eta+1} \end{array}\right](I_{(2\eta+1)\ell}+N) \\ = &  \left[\begin{array}{cccccccc} C_{11} & C_{01} & C_{12} & ... & C_{0\eta} & C_{1(\eta+1)} \end{array}\right],
\end{split}
\end{equation*} 
where
$$ I_{(2\eta+1)\ell}+N:=  \left[\begin{array}{ccccccc} I_\ell & \widehat A \\ & I_\ell  & I_\ell \\ & & I_\ell & \widehat A \\ & & & \ddots & \ddots \\  & & & &  I_\ell & I_\ell \\ & & & & & I_\ell \end{array}\right],$$
and $\widetilde C := C_{0(\eta+1)}-E_{\eta+1}\widehat A$. Clearly
$$ \| \left[\begin{array}{cccccccc} E_1 & F_1 & E_2  & \ldots & F_\eta & E_{\eta+1} \end{array}\right] \|_F =  \|(X_{12},Y_{32})\|_F,$$
$$  \|  \left[\begin{array}{cccccccc} C_{11} & C_{01} & C_{12} & \ldots & C_{0\eta} & C_{1(\eta+1)} \end{array}\right] \|_F \le \|\Delta_{12}(\la)\|_F,
$$
and,  since the matrix $N$ is nilpotent with $N^{2\eta+1}=0$,
$$  (I_{(2\eta+1)\ell}+N)^{-1} = \sum_{i=0}^{2\eta}(-N)^i.
$$
In addition, $N$ has even powers $N^{2i}$ of 2-norm $\|\widehat A^i \|_2\leq \|\widehat A\|_2^i$, whereas the odd powers $N^{2i-1}$ have 2-norm
$\max(\|\widehat A^{i-1}\|_2, \|\widehat A^i\|_2) \leq \max(\|\widehat A\|_2^{i-1}, \|\widehat A\|_2^i)$.  Since both of them can be bounded by $\max(1,\|\widehat A\|_2^{i})$, it then follows that
\begin{equation} \label{ineq4}
\begin{split}
\|(X_{12},Y_{32})\|_F & \le  \|\Delta_{12}(\la)\|_F (1+ 2\,  \max(1,\|\widehat A\|_2)+ \cdots + 2\,  \max(1,\|\widehat A\|_2^{\eta}) )\\
& \le  [1+2\eta \max(1,\|\widehat A\|_2^{\eta})] \, \|\Delta_{12}(\la)\|_F.
\end{split}
\end{equation}

The discussion for the $\widehat B(\la)$ block is clearly analogous and will yield the bound
\begin{equation} \label{ineq5}
\|(X_{23},Y_{21})\|_F  \le [1+2\e \max(1,\|\widehat A\|_2^{\e})] \|\Delta_{21}(\la)\|_F.
\end{equation}

We can thus summarize this discussion in the following Theorem.
\begin{theorem} \label{th:3}
	Let the pencil $\widehat S_2(\la)$ have the anti-triangular form given in \eqref{eqstep3}. Then the updating strict equivalence transformation $(I-X)\widehat S_2(\la)(I-Y)$
	detailed in \eqref{trstep3} exists and can be bounded by
	$$ \|(X,Y)\|_F \le  [1+2 \, \max(\eta,\e) \, \max(1,\|\widehat A\|_2^{\max(\eta,\e)})]\, \| \Delta^{new}_2(\la)\|_F . $$
\end{theorem}
\begin{proof}
	The bound for $\|(X,Y)\|_F$ follows directly from the identity
	$$\|(X,Y)\|_F^2=\|(X_{12},Y_{32})\|_F^2+\|(X_{23},Y_{21})\|_F^2,$$
	from the inequality
	$\|\Delta_{12}(\la)\|_F^2+\|\Delta_{21}(\la)\|_F^2 \le \| \Delta^{new}_2(\la)\|_F^2 $ and from the individual inequalities \eqref{ineq4} and \eqref{ineq5}.
\end{proof}

The following first order bound in Corollary \ref{C3} for the norm of the perturbation error $\Delta^{new}_3(\la)$ follows from Lemma \ref{cum_errors}, Theorem \ref{th:3} and Corollaries \ref{C1} and \ref{C2}.

\begin{corollary} \label{C3}  Let us define the scalar $ f_3:=\sqrt{2}\, [1+2\max(\eta,\e)\max(1,\|\widehat A\|_2^{\max(\eta,\e)})]$. Then	
\begin{equation*}
\begin{split}
	\|\Delta^{new}_3(\la)\|_F &  \le   [1 + f_3 \| \widehat S_2(\la)\|_2] \, \|\Delta_2^{new}(\la)\|_F +{\mathcal O}(\delta^2) \\
 &	\le [1 + f_3 \| S(\la)\|_2]\, \|\Delta_2^{new}(\la)\|_F +{\mathcal O}(\delta^2).
\end{split}	\end{equation*}
\end{corollary}

\subsection{Putting it all together}

In this subsection, we combine the obtained results regarding the strict equivalence transformation that restores in $\widehat{S}(\la)$ of \eqref{eq:perturbed} the special structure of the unperturbed block Kronecker linearization $S(\la)$ defined in \eqref{block_kronecker_lin}, in such a way that the eigenstructure of $\widehat{S}(\la)$ can be linked to that of a particular rational matrix $\widetilde R(\la)$ as in \eqref{nearby_rational}. The final goal is to bound the norms of the differences between the quadruples $\{ \la I_\ell -A, B, C, D(\la)\}$ and $\{ \la I_\ell - \widetilde A, \widetilde B, \widetilde C, \widetilde D(\la)\}$ that are used for representing the unperturbed rational matrix $R(\la)$ and the perturbed one $\widetilde R(\la)$, respectively.

Recall that we were given the pencil $S(\la)$ of which we want to compute the eigenstructure, since it gives the one of the rational matrix $R(\la)$ in \eqref{eq.Rsect4fro}. Instead, our backward stable algorithm applied to $S(\la)$ computes the exact eigenstructure of a slightly perturbed pencil
$\widehat S(\la)$ with additive error $\Delta_{S}(\la)$ which is induced by the
eigenstructure algorithm and is bounded as~:
$$ \|\Delta_{S}(\la)\|_F \le  c(\ell,m\eta,n\e) \cdot \epsilon_M  \cdot \|S(\la)\|_F,
$$
where $\epsilon_M$ is the machine precision of the used computer, and $c(\ell,m\eta,n\e)$ is a moderate function depending only on the size of the matrix pencil.  We then constructed in three steps a new modified block Kronecker linearization
\begin{equation}\label{transformation}
\widetilde S(\la) := (I-X)\widehat S(\la) (I-Y):= (I-X_3)(I-X_2)(I-X_1)\widehat S(\la) (I-Y_1)(I-Y_2)(I-Y_3)
\end{equation}
as in \eqref{eqstep33}, strictly equivalent to $\widehat{S} (\la)$,
where both $\|X\|_F$ and $\|Y\|_F$ are also of the order of the machine precision times some factors and such that the corresponding rational matrix $\widetilde R(\la)$ \eqref{nearby_rational} has a similar representation as $R(\la)$. Since $\widehat S(\la)$ and $\widetilde S(\la)$ are strictly equivalent pencils, they have {\em exactly} the same eigenstructure, which implies that we have
computed the exact eigenstructure of the nearby rational matrix $\widetilde R(\la)$.

For convenience, the blocks of $ \widetilde S(\la) $ will be expressed in the sequel as
$\widetilde{M}(\la):=M(\la)+\Delta M(\la)$, $\widetilde{A}:= A+\Delta A$, $\widetilde{B}:= B+\Delta B$ and $\widetilde{C}:= C+\Delta C$. In the previous subsections, we rewrote $ \widetilde S(\la) $ as an additive perturbation
\[
\widetilde S(\la) = S(\la) + \Delta_3^{new}(\la)
\]
and derived a first order bound for the norm of the error pencil $\Delta_3^{new}(\la)$ in Corollaries \ref{C1}, \ref{C2} and \ref{C3}~:
\begin{equation} \label{fulli}   \|\Delta_3^{new}(\la)\|_F\le (1+f_1\|S(\la)\|_2)(1+f_2\| S(\la)\|_2)(1+f_3\| S(\la)\|_2)\|\Delta_S(\la)\|_F + {\mathcal O}(\delta^2).
\end{equation}
This implies, in particular, that if $\|\Delta_S(\la)\|_F$ is sufficiently small, then the norms of the perturbations $\Delta A$, $\Delta B$ and $\Delta C$ are sufficiently small to guarantee that $\widetilde C(\la I_\ell-	\widetilde A)^{-1}	\widetilde B$ is a minimal state-space realization, as announced. Then, according to \cite{ADMZ}, $ \widetilde S(\la)$ is indeed a strong linearization of the rational matrix $\widetilde{R} (\la)$ in \eqref{nearby_rational}.
Moreover, \eqref{fulli} also implies that if $\|\Delta_S(\la)\|_F$ is sufficiently small, then $ \widetilde D(\la):=\sum_{i=0}^d (D_i+\Delta D_i )\la^i$ in \eqref{nearby_rational} is a polynomial matrix with the same degree $d=\eta+\e+1$ as the polynomial part $D(\la)$ of $R(\la)$ (recall that we are assuming that $d$ is the degree of $D(\la)$ or, equivalently, that $D_d \ne 0$).

Notice that $\widetilde R(\la)$ in \eqref{nearby_rational} is the transfer function of the following perturbed polynomial system matrix
\begin{equation}\label{eq:perturbed_pulsystem}
P(\la)+\Delta P(\la) := \begin{bmatrix}
\la I_\ell-A & -B\\
C & D(\la)
\end{bmatrix}+ \begin{bmatrix}
-\Delta A & -\Delta B\\
\Delta C & \sum_{i=0}^d \Delta D_i \la^i
\end{bmatrix},
\end{equation}
where $P(\la)$ is a polynomial system matrix of the original rational matrix $R(\la)$. Recall that $\|R(\la)\|_F$ is defined in \eqref{def:norm} as $ \|P(\la)\|_F$. This motivates us to define the norm of the perturbation of $R(\la)$ as
$$\|\Delta R(\la)\|_F := \|\Delta P(\la)\|_F = \sqrt{ \| \Delta  A \|_F^2 + \| \Delta  B \|_F^2+ \| \Delta  C \|_F^2+ \sum_{i=0}^d\| {\Delta  D_i}\|_F^2}.$$

After this discussion, we present our main perturbation results in Theorems \ref{th:4} and \ref{th:5}. The first one focuses on block Kronecker linearizations and the second one on the corresponding rational matrices.
\begin{theorem} \label{th:4} Let $R(\la)$ be the $m \times n$ rational matrix in \eqref{eq_rationalmatrix} and let $S(\la)$ be a block Kronecker linearization of $R(\la)$ as in \eqref{block_kronecker_lin}. Let us define $\alpha:=1+2\e \max(1,\| A\|_2^{\e}),$ $\beta:=1+2\eta \max(1,\| A\|_2^{\eta}) ,$ $\gamma:=\frac{\e +\eta}{2\sqrt{2}}$ and $s:= \max(\alpha,\beta, \gamma)+\gamma(\beta\|B\|_2+\alpha\|C\|_2).$  Assume that $\max ( \e , \eta) > 0$ and consider the functions dependent on the initial data
	\begin{equation*}
	\begin{split}
	f_1&:=f_1(\e,\eta,\|A\|_2, \|B\|_2,\|C\|_2):= \frac{4\sqrt{2}s}{2-\sqrt{3}} ,\\
	f_2&:=f_2(\e,\eta) :=\frac{\sqrt{2} \, ( 4\max(\e,\eta)-1)}{3}  ,\\
	f_3&:=f_3(\e,\eta, \| A\|_2) :=\sqrt{2} \, [1+2 \max(\eta,\e) \max(1,\| A\|_2^{\max(\eta,\e)})] .
	\end{split}
	\end{equation*}
	Let $\widehat S(\la):= S(\la) + \Delta_S(\la)$ be a perturbed pencil as in \eqref{eq:perturbed}. If $\|\Delta_{S}(\la)\|_F$ is sufficiently small,
	then $\widehat S(\la)$ is strictly equivalent to a block Kronecker linearization $\widetilde S (\la)$ as in \eqref{eqstep33} with the same parameters $\e$ and $\eta$ as $S(\la)$, i.e., the transformation  \eqref{transformation} exists. Moreover, $\widetilde S(\la) = S(\la) + \Delta_3^{new}(\la)$ with
	\begin{equation} \label{full}  \|\Delta_3^{new}(\la)\|_F \le (1+f_1\|S(\la)\|_2)(1+f_2\|S(\la)\|_2)(1+f_3\|S(\la)\|_2)\|\Delta_S(\la)\|_F + {\mathcal O}(\delta^2),
	\end{equation}
	where $\delta :=\frac{\|\Delta_S(\la)\|_F}{\|S(\la)\|_F}$.
\end{theorem}
\begin{proof}  This follows directly from  \eqref{fulli}, except that we have replaced the 2-norm of $\wh A$ in $f_3$ in Corollary \ref{C3} by that of $A$, because the difference can be absorbed in the ${\mathcal O}(\delta^2)$ term.
\end{proof}


Theorem \ref{th:4} does not provide directly bounds on the norms of the differences between the quadruples representing the rational matrices $R(\la)$ and $\widetilde{R}(\la)$ corresponding to the block Kronecker linearizations $S(\la)$ and $\widetilde{S} (\la)$. The reason is that the polynomial parts  $D(\la) = (\Lambda_\eta(\la)\otimes I_m)^T M(\la) (\Lambda_\epsilon(\la)\otimes I_n)$ and  $\widetilde D(\la) = (\Lambda_\eta(\la)\otimes I_m)^T \widetilde M(\la) (\Lambda_\epsilon(\la)\otimes I_n)$ of $R(\la)$ and $\widetilde{R}(\la)$ are not directly visible in $S(\la)$ and $\widetilde{S}(\la)$. For this reason, we will need Lemma \ref{le:D}, that follows from \cite[Lemma 2.15, Theorem 4.4 and Lemma 5.23(b)]{DLPV}.
\begin{lemma} \label{le:D}
	Let $M(\la)$ be a $m(\eta +1) \times n(\e +1)$ pencil and let $\Lambda_k(\la) := \left[\begin{array}{ccccc} \la^k & \cdots & \la & 1 \end{array}\right]^T.$
	If we define the polynomial matrix $Q(\la)$ as
	\begin{equation} \label{DM} Q(\la) :=(\Lambda_\eta(\la)\otimes I_m)^T \, M (\la) \, (\Lambda_\e(\la)\otimes I_n),
	\end{equation}
	then we can bound its norm as follows
	$$  \| Q(\la)\|_F \le \sqrt{2\min(\e+1,\eta+1)} \, \|M(\la)\|_F.
	$$
	Moreover, for every polynomial matrix $Q(\la)$ of degree at most $d = \e + \eta + 1$, there exist infinitely many pencils $M (\la)$ satisfying \eqref{DM}. For each of these pencils $\|M(\la) \|_F \geq \| Q(\la)\|_F / \sqrt{2 d}$ and there exist pencils such that $\| Q(\la)\|_F =  \|M (\la)\|_F$.
\end{lemma}
As commented in \cite{DLPV}, Fiedler and proper generalized Fiedler pencils (modulo permutations) of a polynomial matrix $Q(\la)$ satisfy $\| Q(\la)\|_F =  \|M (\la)\|_F$ in Lemma \ref{le:D}. On the other hand, it might be worth to remind that there exist pencils $M(\la)$ satisfying \eqref{DM} with norm arbitrarily larger than the norm of $Q(\la)$.


We are finally in the position of proving the main perturbation result of this paper.	
\begin{theorem}\label{th:5} Let $R(\la) = C(\la I_\ell-A)^{-1}B + \sum_{i=0}^d D_i\la^i$ be an $m \times n$ rational matrix, where $C(\la I_\ell-A)^{-1}B$ is a minimal state-space realization of the strictly proper part of $R(\la)$, let $S(\la)$ be a block Kronecker linearization of $R(\la)$ as in \eqref{block_kronecker_lin} with $\max ( \e , \eta) > 0$, and let $f_1, f_2, f_3$ be the functions defined in Theorem \ref{th:4}. Let $\widehat S(\la):= S(\la) + \Delta_S(\la)$ be a perturbed pencil as in \eqref{eq:perturbed}. If $\|\Delta_{S}(\la)\|_F$ is sufficiently small,
	then $\widehat S(\la)$ is strictly equivalent to a block Kronecker linearization $\widetilde S (\la)$ as in \eqref{eqstep33},  with the same parameters $\e$ and $\eta$ as $S(\la)$, of a rational matrix
	$$
	\widetilde R(\la) = \widetilde  C(\la I_\ell- \widetilde  A)^{-1} \widetilde B + \sum_{i=0}^d \widetilde D_i\la^i,
	$$
	where $\widetilde C(\la I_\ell-\widetilde A)^{-1} \widetilde B$ is a minimal state-space realization of the strictly proper part of $\widetilde R(\la)$. Moreover, if $\widetilde{A}:= A+\Delta A$, $\widetilde{B}:= B+\Delta B$,  $\widetilde{C}:= C+\Delta C$ and $\widetilde{D}_i := D_i +\Delta D_i$, $i = 0, 1, \ldots , d$, then		
	\begin{equation} \label{eq.KSRbound}
	\dfrac{\sqrt{ \| \Delta  A \|_F^2 + \| \Delta  B \|_F^2+ \| \Delta  C \|_F^2+ \sum_{i=0}^d\| {\Delta  D_i}\|_F^2}}{\|R(\la)\|_F}\le K_{S,R} \dfrac{\|\Delta_S(\la)\|_F}{\|S(\la)\|_F} + {\mathcal O}( \delta^2 ) ,
	\end{equation}
	where $$K_{S,R}:= \sqrt{2\min(\e+1,\eta+1)}\, (1+f_1\|S(\la)\|_2)(1+f_2\|S(\la)\|_2)(1+f_3\|S(\la)\|_2)\frac{\|S(\la)\|_F}{\|R(\la)\|_F}$$
	and $\delta = \frac{\|\Delta_S(\la)\|_F}{\|S(\la)\|_F}$.
\end{theorem}
\begin{proof}
	Since $\widetilde{S} (\la)$ and $S(\la)$ have the same structure according to Theorem \ref{th:4},
	\[
	\Delta_3^{new} (\la) = \widetilde{S} (\la) - S(\la) =
	\left[\begin{array}{ccc} \widetilde M(\la) - M(\la) & \widehat K_2^T (\widetilde C -C) &  0 \\ (\widetilde B - B) \widehat K_1 & \widetilde A- A & 0\\
	0 &   0 &  0
	\end{array}\right]
	\]
	and $\|\Delta_3^{new} (\la)\|_F = \sqrt{\| \Delta  A \|_F^2 + \| \Delta  B \|_F^2+ \| \Delta  C \|_F^2+ \|  \widetilde M(\la) - M(\la) \|_F^2}$. Next, we combine this expression of $\|\Delta_3^{new} (\la)\|_F$ with
	$\sum_{i=0}^d D_i\la^i = (\Lambda_\eta(\la)\otimes I_m)^T M(\la) (\Lambda_\epsilon(\la)\otimes I_n)$,  $\sum_{i=0}^d \widetilde D_i\la^i = (\Lambda_\eta(\la)\otimes I_m)^T \widetilde M(\la) (\Lambda_\epsilon(\la)\otimes I_n)$ and Lemma \ref{le:D}, and we get
	\[
	\sqrt{ \| \Delta  A \|_F^2 + \| \Delta  B \|_F^2+ \| \Delta  C \|_F^2+ \sum_{i=0}^d\| {\Delta  D_i}\|_F^2} \leq \sqrt{2\min(\e+1,\eta+1)}\, \|\Delta_3^{new} (\la)\|_F \, .
	\]
	The rest of the proof follows from \eqref{full}.
\end{proof}

The strength of the new structured backward error analysis that we present in this paper for the computation of the eigenstructure of a rational matrix $R(\la)$ by applying a backward stable generalized eigenvalue algorithm to a block Kronecker linearization $S(\la)$ of $R(\la)$ is that we can interpret the computed eigenstructure as the exact eigenstructure for a slightly perturbed rational matrix $\widetilde R(\la)$ corresponding to the nearby quadruple $\{\la I_\ell- \widetilde A, \widetilde B, \widetilde C, \widetilde D(\la) \},$
and that we have a bound on the error because we have a specific coordinate system in which we can describe both the original rational matrix $R(\la)$ and its perturbed version $\widetilde R(\la)$, namely by the quadruples $\{\la I_\ell-A,B,C,D(\la)\}$ and $\{\la I_\ell-\widetilde A,\widetilde B,\widetilde C,\widetilde D(\la)\}$. It still remains to analyze under which conditions this bound is satisfactory. This is the purpose of the next subsection.

\subsection{Sufficient conditions for structural backward stability} \label{subsec.sufficientconds}
The goal of this section is to establish sufficient conditions on $R(\la)$ and $S(\la)$ that guarantee that $K_{S,R}$ in \eqref{eq.KSRbound} is moderate and, thus, that guarantee structural backward stability. We advance that these conditions are the following
\begin{equation}\label{eq.sufficient_conditions}
\max( \|A\|_F , \|B\|_F , \|C\|_F, \|D(\la)\|_F ) \leq 1 \quad \mbox{and} \quad
\|M(\la)\|_F \approx \|D(\la)\|_F,
\end{equation}
where the notation introduced in the previous section is used. Observe that the first condition is a condition on $R(\la)$ while the second one is on $S(\la)$. According to Lemma \ref{le:D}, the second condition can be satisfy simply by choosing an adequate block Kronecker linearization $S(\la)$. In addition, we will see that the conditions \eqref{eq.sufficient_conditions} are essentially necessary for $K_{S,R}$ to be moderate, though this does not mean that they are necessary for structural backward stability since \eqref{eq.KSRbound} is an upper bound. For the sake of clarity, the discussion in this section focuses on identifying the key ingredients for structural backward stability instead of on providing precise bounds. There exist, obviously, rational matrices which do not satisfy the first condition in \eqref{eq.sufficient_conditions}. We will discuss in Section \ref{scaling} how to proceed in such cases.

In the first place observe that each of the essential four factors of $K_{S,R}$, that is, $(1+f_1\|S(\la)\|_2)$, $(1+f_2\|S(\la)\|_2)$, $(1+f_3\|S(\la)\|_2)$ and $\frac{\|S(\la)\|_F}{\|R(\la)\|_F}$, is larger than $1$.
This is obvious for the first three factors. For the fourth factor, it follows from the equalities
\begin{equation} \label{eq.normSnormR}
\begin{split} \|S(\la)\|_F^2 & =   \|A\|_F^2  + \|B\|_F^2  + \|C\|_F^2  + \|M(\la)\|_F^2  + \ell + 2(m\eta+n\e) \quad \text{and} \\
\|R(\la)\|_F^2 & =    \|A\|_F^2  + \|B\|_F^2  + \|C\|_F^2  +  \|D(\la)\|_F^2 + \ell \, .
\end{split}
\end{equation}

To find upper bounds for the three factors $(1+f_1\|S(\la)\|_2)$, $(1+f_2\|S(\la)\|_2)$, $(1+f_3\|S(\la)\|_2)$ of $K_{S,R}$ requires to upper bound each $f_i$ and $\|S(\la)\|_2$. For this purpose, we consider Lemmas \ref{bounds} and \ref{lemm.boundSla}. Lemma \ref{bounds} provides a bound on the function $f_1$ that allows us to identify its most relevant dependencies. Moreover, Lemma \ref{bounds} emphasizes the key role of $t:=\max(\eta,\e)$ in our perturbation analysis. Lemma \ref{lemm.boundSla} bounds $\|S(\la)\|_2$.
\begin{lemma}\label{bounds}
	Let us define $ M_a:= \max(1,\|A\|_2)$,  $M_b:= \max(\|B\|_2,\|C\|_2)$ and $t:=\max(\eta,\e) > 0$ and consider the functions $f_1$, $f_2$ and $f_3$ in Theorem \ref{th:4}. Then
	$$1 \leq f_1 \le 22(1+2\, t M_a^t)(1+ \sqrt{2} \, t M_b), \,\,\, 1 \leq f_2 = \dfrac{\sqrt{2}}{3}(4t-1), \,\,\, 1 \leq f_3 = \sqrt{2} \, (1+2\, t M_a^t).
	$$
\end{lemma}
\begin{proof} It follows by taking into account the inequalities
	$ \gamma \le \frac{t}{\sqrt{2}}$ and $ s \le (1+2 t M_a^t)(1+ \sqrt{2} \, t \, M_b).$
\end{proof}

\begin{lemma} \label{lemm.boundSla} Let $S(\la)$ be the block Kronecker linearization \eqref{block_kronecker_lin}. Then $$\max(1, \|A\|_2 , \|B\|_2 , \|C\|_2, \|M(\la)\|_2 ) \leq \|S(\la)\|_2 $$ and 
{\small$$\|S(\la)\|_2 \leq \sqrt{2} + \|
	\left[\begin{array}{cc} M(\la) & \widehat K_2^T C   \\ B \widehat K_1 & A  \\
	\end{array}\right] \|_2 \leq \sqrt{2} +
	\sqrt{\|A\|_F^2  + \|B\|_F^2  + \|C\|_F^2  + \|M(\la)\|_F^2}  .$$}
\end{lemma}

\begin{proof} The first inequality follows from the definition of the 2-norm of a pencil given in the introduction and the fact that the 2-norm of a matrix is larger than or equal to the 2-norm of any of its submatrices. The second inequality follows from applying the triangular inequality to
	$$
	S(\la) = \left[\begin{array}{ccc} M(\la) & \widehat K_2^T C & 0 \\ B \widehat K_1 & A  & 0\\
	0 &   0 &  0
	\end{array}\right]+
	\left[\begin{array}{ccc} 0 & 0 &  K_2^T(\la)  \\ 0 &  - \la I_\ell & 0\\
	K_1(\la) &   0 &  0
	\end{array}\right] \, .
	$$
	Note that the 2-norm of a pencil as defined in the introduction is indeed a norm and, so, the triangular inequality can be applied.
\end{proof}

We remark that Lemmas \ref{bounds} and \ref{lemm.boundSla} imply that the conditions \eqref{eq.sufficient_conditions} are essentially necessary for $K_{S,R}$ to be moderate. This can be seen as follows. First, from Lemma \ref{le:D}, we have $\|M(\la) \|_F \geq \|D(\la)\|_F / \sqrt{2(\e +\eta +1)}$. Thus, $\max( \|A\|_F , \|B\|_F , \|C\|_F, \allowbreak \|D(\la)\|_F ) \gg 1$ implies $\|S(\la)\|_2 \gg 1$, which in turns implies $K_{S,R} \gg 1$, since $f_i \geq 1$ for $i=1,2,3$. Moreover, if $\|M(\la)\|_F \gg \|D(\la)\|_F$, then $\|S(\la)\|_F / \|R(\la)\|_F \gg 1$ may happen, according to \eqref{eq.normSnormR}, and $K_{S,R} \gg 1$ in that situation. We emphasize that the condition $\|M(\la)\|_F \approx \|D(\la)\|_F$ was also used in the analysis in \cite[Corollary 5.24]{DLPV}.

Next, we prove the announced result that conditions \eqref{eq.sufficient_conditions} are sufficient for $K_{S,R}$ to be moderate and, thus, for structural backward stability.

\begin{corollary} \label{cor.definitive} Under the hypotheses and with the notation of Theorem \ref{th:5}, assume, in addition, that \eqref{eq.sufficient_conditions} holds and let $t:=\max(\eta,\e) > 0$. Then,
	$$
	K_{S,R} \leq \, g \, t^q \sqrt{m+n} \, ,
	$$
	where $q =5$, if $\eta > 0$ and $\e >0$, $q = 9/2$, if $\eta = 0$ or $\e = 0$, and $g$ is a moderate number (a constant that does not depend on $\eta, \e, m, n, \ell$). Moreover
	$$
	\dfrac{\sqrt{ \| \Delta  A \|_F^2 + \| \Delta  B \|_F^2+ \| \Delta  C \|_F^2+ \sum_{i=0}^d\| {\Delta  D_i}\|_F^2}}{\|R(\la)\|_F}\leq \, g \, t^q \sqrt{m+n} \,\, \dfrac{\|\Delta_S(\la)\|_F}{\|S(\la)\|_F} + {\mathcal O}( \delta^2 ) \, .
	$$
\end{corollary}
\begin{proof}
	Note that \eqref{eq.sufficient_conditions} and Lemmas \ref{bounds} and \ref{lemm.boundSla} imply $\|S(\la)\|_2 \lesssim 2 + \sqrt{2}$, $f_1 \leq g_1 t^2$, $f_2 \leq g_2 t$, and $f_3 \leq g_3 t$, with $g_1, g_2, g_3$ moderate numbers. Moreover, from \eqref{eq.normSnormR}, \eqref{eq.sufficient_conditions} and $\|R(\la)\|_F \geq 1$, we get that
	$\|S(\la)\|_F^2 \approx \|R(\la)\|_F^2   + 2(m\eta+n\e)$ and
	$$
	\|S(\la)\|_F^2 \leq (1+2(m\eta+n\e))\, \|R(\la)\|_F^2 \leq \, 3 \, (m+n) \, t \, \|R(\la)\|_F^2.
	$$
	It only remains to analyze the factor $\sqrt{2\min(\e+1,\eta+1)}$ of $K_{S,R}$, which is less than or equal to $\sqrt{2 (t+1)}$, if $\eta > 0$ and $\e > 0$, or equal to $\sqrt{2}$, if
	$\eta = 0$ or $\e = 0$. Combining all these bounds with the fact that $t \geq 1$, the result follows as a corollary of Theorem \ref{th:5}.
\end{proof}

\begin{remark} Observe that \eqref{eq.sufficient_conditions} allow $\max( \|A\|_F , \|B\|_F , \|C\|_F, \|D(\la)\|_F )  \ll 1$. However, since the rational matrix $R(\la)$ in \eqref{eq_rationalmatrix} can be multiplied by a nonzero number without affecting at all its eigenstructure, it is natural and convenient to use as sufficient conditions
	\begin{equation}\label{eq.sufficient_conditions2}
	\max( \|A\|_F , \|B\|_F , \|C\|_F, \|D(\la)\|_F ) = 1 \quad \mbox{and} \quad
	\|M(\la)\|_F \approx \|D(\la)\|_F.
	\end{equation}
	Such conditions would have appeared as sufficient in the analysis if we had defined the norm of $R(\la)$ as
	\begin{equation} \label{eq.secondnormofR}
	|||R(\la)|||_F := \sqrt{\| A \|_F^2 + \| B \|_F^2+ \| C \|_F^2+ \sum_{i=0}^d\| {D_i}\|_F^2},
	\end{equation}
	instead as in \eqref{def:norm} (observe that we have removed the $\ell$ summand), depending only on the free parameters of the representation of $R(\la)$ in \eqref{eq_rationalmatrix}. We have chosen to use \eqref{def:norm} because, first, it identifies the informal ``norm'' of $R(\la)$ with the formal norm of the polynomial system matrix $P(\la)$ and, second, it corresponds to the particular case $E = I_\ell$ of the more general representation $R(\la)=C(\la E - A)^{-1}B + D(\la),$ with $E$ nonsingular, when taking as norm the one of the corresponding polynomial system matrix. Under the conditions \eqref{eq.sufficient_conditions2}, it is essentially equivalent to use \eqref{def:norm} or \eqref{eq.secondnormofR} as ``norm'' of $R(\la)$. The use of representations $R(\la)=C(\la E - A)^{-1}B + D(\la)$ for rational matrices is of interest in certain applications and the block Kronecker linearizations in this case are obtained just by replacing $A - \la I_\ell$ by $A - \la E$ in \eqref{block_kronecker_lin}. We will consider the analysis of this general case in the future.
\end{remark}

\subsection{Restoring the structure when the polynomial part of the rational matrix is linear} \label{subsec.linearD}

In this subsection, we consider the particular case of having a rational matrix with linear polynomial part. That is, the case of having a rational matrix that can be written in the form
$$R(\la)=C(\la I_\ell - A)^{-1}B + M(\la),$$
where $C(\la I_\ell - A)^{-1}B$ is a minimal state-space realization and $M(\la)$ is a matrix pencil. Then $R(\la)$ can be strongly linearized using the following linear polynomial system matrix
\begin{equation}\label{system}
S(\la):=\left[\begin{array}{ccc} M(\la) & C \\   B & A- \la I_\ell
\end{array}\right].
\end{equation}
Notice that, in this case, the linearization does not have the block anti-triangular structure as the block Kronecker linearization in \eqref{block_kronecker_lin} since $K_{1}(\la)$ and $K_{2}(\la)$ are empty matrices. The strong linearization \eqref{system} can be seen as the limit case of \eqref{block_kronecker_lin} when $\e = \eta = 0$.

If we compute the eigenstructure of $S(\la)$, the backward stability of the staircase algorithm \cite{VanDooren79} and the $QZ$ algorithm \cite{MOS} guarantees that we computed the exact eigenstructure of a slightly perturbed pencil
\begin{equation}\label{eq:perturbed2}
\widehat S(\la):= S(\la) + \Delta_S(\la), \quad  \Delta_S(\la) := \left[\begin{array}{cc} \Delta_{11}(\la) & \Delta_{12}(\la)  \\ \Delta_{21}(\la) & \Delta_{22}(\la)
\end{array}\right].
\end{equation}
The structure of \eqref{system} is lost in \eqref{eq:perturbed2} since the off-diagonal blocks of $\widehat S(\la)$ are not constant matrices and the identity block $I_\ell$  is not preserved by the perturbation.

Notice that restoring in $\widehat S(\la)$ the original structure of $S(\la)$ is much simpler than in previous sections, as we do not have to restore any anti-triangular zero block nor the minimal bases $K_{1}(\la)$ and $K_{2}(\la)$ in \eqref{eq:pencil}. We only have to take care of restoring the identity matrix $I_\ell$ and the constant matrices $B$ and $C$ to obtain in two steps a new strictly equivalent linear polynomial system matrix
\begin{equation}\label{transformation2}
\widetilde S(\la) := (I-X)\widehat S(\la) (I-Y):= (I-X_2)(I-X_1)\widehat S(\la) (I-Y_1)(I-Y_2)
\end{equation} of the form

\begin{equation}\label{systemtilde}
\widetilde S(\la):=\left[\begin{array}{ccc} \widetilde  M(\la) & \widetilde  C \\   \widetilde  B & \widetilde  A- \la I_\ell
\end{array}\right],
\end{equation}
where $\widetilde{M}(\la):=M(\la)+\Delta M(\la)$, $\widetilde{A}:= A+\Delta A$, $\widetilde{B}:= B+\Delta B$ and $\widetilde{C}:= C+\Delta C$.
For that, we consider the discussion in Subsection \ref{restoring_2}, for restoring $I_\ell$; and a simplified version of the discussion in Subsection \ref{restoring_3}, for restoring the constant matrices $B$ and $C$. In particular, from the bound in \eqref{ineq3} and a counterpart of Theorem \ref{th:3} we get the following result.

\begin{theorem} \label{th:5linearpolynomialpart}
	
	Let $S(\la)$ be a minimal linear system matrix as in \eqref{system}. The transformation $(X,Y)$ in \eqref{transformation2} exists and we can bound the corresponding perturbation $\widetilde S(\la) - S(\la)$
	as follows~:
	\begin{equation} \label{full2}  \| \widetilde S(\la) - S(\la) \|_F \le (1+\sqrt{2}\|S(\la)\|_2)^2 \, \|\Delta_S(\la)\|_F + {\mathcal O}(\delta^2).
	\end{equation}
	In addition, if $\| \widetilde S(\la) - S(\la) \|_F$ is sufficiently small, then the perturbed pencil $\widetilde S(\la)$ is a minimal linear system matrix of the rational matrix $\widetilde R(\la) = \widetilde C(\la I_\ell-	\widetilde A)^{-1}	\widetilde B + \widetilde M(\la)$ and
	\[
	\dfrac{\sqrt{ \| \Delta  A \|_F^2 + \| \Delta  B \|_F^2+ \| \Delta  C \|_F^2+ \| \Delta M(\la) \|_F^2}}{\|R(\la)\|_F}\le  (1+\sqrt{2}\|S(\la)\|_2)^2 \, \dfrac{\|\Delta_S(\la)\|_F}{\|S(\la)\|_F} + {\mathcal O}(\delta^2) ,
	\]
	where $\delta = \|\Delta_S(\la)\|_F / \|S(\la)\|_F$.
\end{theorem}
The simplicity of the bound in Theorem \ref{th:5linearpolynomialpart} is also a consequence of $\|S(\la)\|_F = \|R(\la)\|_F$.

\section{Scaling for obtaining structural backward stability} \label{scaling}
Once a block Kronecker linearization $S(\la)$ in \eqref{block_kronecker_lin} of $R(\la)$ in \eqref{eq_rationalmatrix} satisfying $\|M(\la)\|_F \approx \|D(\la)\|_F$ is chosen and the staircase or the $QZ$ algorithm is applied to $S(\la)$, structural backward stability is guaranteed for the computed eigenstructure if the first condition in \eqref{eq.sufficient_conditions} holds. However, there exist rational matrices which do not satisfy $\max( \|A\|_F , \|B\|_F , \|C\|_F, \|D(\la)\|_F ) \leq 1$ and, therefore, the computation of their eigenstructure via a block Kronecker linearization might not be structurally backward stable. In this section, we study how to proceed in these cases.


First observe that the eigenstructure of the rational matrix $R(\la)$ does not change at all if it is multiplied by a positive real constant $d_R$. Choosing appropriately $d_R$, we get easily a rational matrix such that $\max( \|B\|_F , \|C\|_F, \|D(\la)\|_F ) \leq 1$. Even more, if $d_R$ is an integer power of $2$, this multiplication can be performed without introducing any rounding error. This indicates that the crucial point is how to deal with rational matrices with $\|A\|_F > 1$. For this, note that when representing a rational matrix $R(\la)$ by a realization quadruple $\{\la I_\ell- A,B,C,D(\la)\}$, where $D(\la)$ is polynomial,
$$ R(\la) := C(\la I_\ell-A)^{-1}B + \sum_{i=0}^d D_i\la^i,
$$
one can change the coordinate system of the state-space realization $\{A,B,C\}$ of the strictly proper part of $R(\la)$ by a diagonal similarity scaling $T:= \diag{d_1, \ldots, d_\ell}$, $d_i > 0$, without changing $R(\la)$ since
$$ C(\la I_\ell-A)^{-1}B = CT(\la I_\ell-T^{-1}AT)^{-1}T^{-1}B .
$$
Thus, before multiplying $R(\la)$ by $d_R$, we can choose $T$ to balance $A$, i.e., to minimize its Frobenius norm under all diagonal similarities by making the 2-norms of the rows and columns of $T^{-1}AT$ become equal \cite{Parlett}. Moreover, at the same time, the Frobenius norms of $T^{-1} B$ and $C T$ can be made equal by considering a positive scalar factor multiplying $T$. Observe, in addition, that if the entries of $T$ are integer powers of $2$, this process does not introduce rounding errors, though, in this case, the norm of $T^{-1}AT$ is only approximately minimized. However, the effects of $T$ are limited since $\|T^{-1}AT\|_F \geq \sqrt{|\la_1|^2 + \cdots + |\la_\ell|^2}$, where $\la_1 , \ldots , \la_\ell$ are the eigenvalues of $A$, for any invertible $T$, i.e., diagonal or not. Therefore, other approaches are needed for dealing with all instances of matrices $A$ with large norms. It is important to emphasize at this point that the influence of a large norm matrix $A$ on the bound \eqref{eq.KSRbound} is huge, because it contributes to $\|S(\la)\|_2$, but also the factor $\|A\|_2^{\max(\eta,\e)}$ is present in both $f_1$ and $f_3$.

The final solution comes from changing the variable $\la$ to $\widehat \la := d_\la \la$ and from combining this with the multiplication by the constant $d_R$ and the diagonal scaling $T$ discussed above. Note that the change of variable transforms the zeros and the poles of $R(\la)$ in a very simple way, preserving their partial multiplicities, and that does not change at all its minimal indices \cite{MMMM-mobius,VanDoorenPhDThesis}. The combination of all these scalings  yields a new transfer function
\begin{equation} \label{eq.scaledhatR}
\widehat R(\widehat \la) := \widehat D(\widehat \la) + \widehat C(\widehat \la I_\ell-\widehat A)^{-1}\widehat B :=   d_R R(\widehat \la/d_\la)
\end{equation}
where
\begin{equation} \label{eq.scaledhatcoefs}
\widehat A := d_\la T^{-1}A T, \;\; \widehat B :=   \sqrt{d_\la d_R}\, T^{-1}B, \;\; \widehat C :=   \sqrt{d_\la d_R}\, C T
\end{equation}
and
\begin{equation} \label{eq.scaledhatcoefs2}
\widehat D_i := d_Rd_\la^{-i} D_i, \quad \text{for all} \quad i=0, 1, \ldots ,d.
\end{equation}
Then, we can choose $d_\la := \min(1,\|T^{-1} A T\|_F^{-1})$, such that $\widehat A$ has norm smaller than or equal to $1$.  Note that the preliminary balancing will make this step milder, in the sense that $d_\la$ will be closer to 1.  Finally, based on \eqref{eq.scaledhatR}, we summarize the following scaling procedure for obtaining a rational matrix $\widehat R(\widehat \la)$ with $\max( \|\widehat A\|_F , \|\widehat B\|_F , \|\widehat C\|_F, \|\widehat D(\widehat \la)\|_F ) = 1$ from the data $\{A, B,C, D_0, D_1, \ldots , D_d\}$:
\begin{description}
	\item[\bf Step 1.] Compute $T = \diag{d_1, \ldots, d_\ell}$ to balance $A$ and to make equal the norms of $T^{-1}B$ and $CT$.
	\item[\bf Step 2.] Choose $d_\la := \min(1,\|T^{-1} A T\|_F^{-1})$.
	\item[\bf Step 3.] Choose
	$$
	d_R = \frac{1}{\displaystyle \max (\,  \| \sqrt{d_\la}\, T^{-1}B \|_F^2 , \|\sqrt{d_\la}\, C T\|_F^2 , \sqrt{\sum_{i=0}^{d}  \|d_\la^{-i} D_i \|_F^2 } \, )}.
	$$
	\item[\bf Step 4.] Compute $\{\widehat A, \widehat B, \widehat C, \widehat D_0, \widehat D_1, \ldots , \widehat D_d\}$ as in \eqref{eq.scaledhatcoefs}-\eqref{eq.scaledhatcoefs2}.
\end{description}
This process can be easily arranged to use scale factors that are all integer powers of two and, thus, can be implemented without any rounding error. Moreover, this scaling can be applied directly to the pencil $S(\la)$. More precisely, the pencil
$$ \widehat S(\widehat \la) := D_\ell S(\widehat \la/d_\la) D_r,
$$
where the left and right diagonal scalings $D_\ell$ and $D_r$ are given by
$$ D_\ell := \diag{d_R^\frac12 d_\la^{-\eta} I_m, \ldots , d_R^\frac12 d_\la^0 I_m, d_\la^\frac12 d_1^{-1}, \ldots , d_\la^\frac12 d_\ell^{-1}, d_R^{-\frac12} d_\la^{\epsilon} I_n,\ldots, d_R^{-\frac12} d_\la^{1} I_n},
$$
$$ D_r := \diag{d_R^\frac12 d_\la^{-\epsilon} I_n, \ldots , d_R^\frac12 d_\la^0 I_n, d_\la^\frac12 d_1, \ldots , d_\la^\frac12 d_{\ell}, d_R^{-\frac12} d_\la^{\eta} I_m,\ldots, d_R^{-\frac12} d_\la^{1}I_m},
$$
is a block Kronecker linearization of the rational matrix $\widehat R (\widehat \la)$ in \eqref{eq.scaledhatR}.




\section{Numerical experiments} \label{numerical}
In this section, we describe three experiments illustrating that the potential sources of
structural backward instability revealed by the bound \eqref{eq.KSRbound} are indeed observed in practice. More precisely, the experiments will illustrate that if a rational matrix $R(\la)$ as in \eqref{eq_rationalmatrix} does not satisfy the first condition in \eqref{eq.sufficient_conditions}, then the computation of the eigenstructure of $R(\la)$ by applying the $QZ$ algorithm to a block Kronecker linearization $S(\la)$ of $R(\la)$ that satisfies $\|M(\la)\|_F = \|D(\la)\|_F$ is not structurally backward stable. Moreover, the experiments also illustrate that
the scaling described in Section \ref{scaling} is effective and leads to structured backward stability for the scaled rational matrices and linearizations.

A difficulty for performing fully reliable numerical experiments in this setting is that to estimate the actual global backward error for the {\em whole} computed eigenstructure, i.e., the left-hand side of \eqref{eq.KSRbound}, is a challenging optimization problem for which we do not know yet a solution.
Therefore, we will limit ourselves to computing a lower bound for the backward error based on the ``local'' backwards errors of each computed zero of the rational matrix, as we explain below. This lower bound might severely underestimate the actual global backward error. Thus, we cannot check from our experiments the sharpness of the bound \eqref{eq.KSRbound}, which, on the other hand, was deduced through many potentially overestimating inequalities with the main goal of getting a bound as clear as possible instead of optimizing its sharpness.

For simplicity, we will restrict our numerical experiments to square and regular rational matrices $R(\la)$ with a corresponding quadruple $\{A,B,C,D(\la)\}$ of moderate dimensions and degree of its polynomial part: $m=n=2$, $\ell = 5$, $d=3$.
The block Kronecker pencil we choose for our computations is
$$   S(\la):=\left[\begin{array}{cccc} \la D_3+D_2 & 0  &  0  &  I_2  \\
0 & \la D_1 + D_0 & C & - \la I_2 \\
0 & B & A -\la I_\ell  & 0\\
I_2 & -\la I_2 &  0 &  0
\end{array}\right], $$
which has $\eta$ and $\e$ equal to $1$, size $11\times 11$ and satisfies $\|M(\la)\|_F = \|D(\la)\|_F$.
We also will look at the polynomial system matrix
$$  P(\la) := \left[ \begin{array}{cc} A-\la I_\ell & B \\ C & D(\la) \end{array}\right], \quad D(\la):=D_0+\la D_1 + \la^2 D_2 + \la^3 D_3
$$
of $R(\la)$ because it allows us to estimate the backward errors of our algorithm as follows. We look for a rational matrix $\widetilde R(\la)$ corresponding to a quadruple $\{A + \Delta A,B + \Delta B,C + \Delta C, (D+ \Delta D)(\la)\}$ such that all its finite zeros are exactly all the computed finite eigenvalues obtained by applying the $QZ$ algorithm to $S(\la)$ and such that $\|(\, \Delta A,\Delta B,\Delta C, (\Delta D)(\la) \, )\|_F$ is as small as possible. As a consequence of the classical results of Rosenbrock \cite{Ros70}, this is equivalent to find a perturbed polynomial system matrix $ P(\la) +\Delta P(\la) $ of $\widetilde R(\la)$, whose finite zeros are the computed eigenvalues $\la_i$ and such that $\|(\, \Delta A,\Delta B,\Delta C, (\Delta D)(\la) \, )\|_F$ is as small as possible. Therefore, $\{ \Delta A, \Delta B,
\Delta C, \Delta D_0, \Delta D_1, \Delta D_2, \Delta D_3 \}$ must have the property that {\em simultaneously}, at each computed eigenvalue $\la_i$, the matrix
{\small $$P(\la_i) +\Delta P(\la_i) = P(\la_i)+  \left[\begin{array}{cc|cc|cc|cc} \Delta A & \Delta B & 0 & 0 & 0 & 0 & 0 & 0 \\
\Delta C & \Delta D_0 & 0 & \Delta D_1 & 0  & \Delta D_2 & 0  & \Delta D_3  \end{array}\right]
\left[\begin{array}{cc} I_{\ell+m} \\  \la_i I_{\ell+m} \\ \la_i^2  I_{\ell+m} \\  \la_i^3  I_{\ell+m} \end{array}\right]$$}
must be singular. To find the smallest possible Frobenius norm of all possible $\{ \Delta A, \Delta B, \allowbreak \Delta C, \Delta D_0, \allowbreak \Delta D_1, \Delta D_2, \Delta D_3 \}$ that satisfy this property {\em for all} computed $\la_i$ is not obvious, however to solve this problem {\em for only one} computed $\la_i$ is easy. For this purpose, let $\Delta^{(i)}$ be the minimum Frobenius norm matrix that makes
$P(\la_i) +\Delta^{(i)}$ singular. Note that $\Delta^{(i)}$ can be computed through the singular value decomposition of $P(\la_i)$ and that, generically, it is a rank one matrix  with Frobenius norm equal to $\sigma_{\min}P(\la_i)$. Then, the linear system
{\small $$ \Delta^{(i)} :=  \left[\begin{array}{cc} \Delta^{(i)}_{11} & \Delta^{(i)}_{12} \\ \Delta^{(i)}_{21} & \Delta^{(i)}_{22} \end{array}\right]
= \left[\begin{array}{cc|cc|cc|cc} \Delta A & \Delta B & 0 & 0 & 0 & 0 & 0 & 0 \\
\Delta C & \Delta D_0 & 0 & \Delta D_1 & 0  & \Delta D_2 & 0  & \Delta D_3  \end{array}\right]
\left[\begin{array}{cc} I_{\ell+m} \\  \la_i I_{\ell+m} \\ \la_i^2  I_{\ell+m} \\  \la_i^3  I_{\ell+m} \end{array}\right]
$$}
for the unknowns
$\{ \Delta A, \Delta B, \Delta C, \Delta D_0, \allowbreak \Delta D_1, \Delta D_2, \Delta D_3 \}$ is consistent and its minimum Frobenius norm solution is given by
$$ \Delta A:=\Delta^{(i)}_{11}, \; \; \Delta B:=\Delta^{(i)}_{12}, \; \; \Delta C:=\Delta^{(i)}_{21}, \; \;  \Delta D_k :=\Delta^{(i)}_{22}\overline \la_i^k/g(\la_i), \; k=0,1,2,3,
$$
where $g(\la_i):=(1+|\la_i|^2+|\la_i|^4+|\la_i|^6)$, and the Frobenius norm of this 7-tuple of matrices is given by
$$ r(P,\la_i) := \| \left[\begin{array}{cc} \Delta^{(i)}_{11} & \Delta^{(i)}_{12} \\ \Delta^{(i)}_{21} & \Delta^{(i)}_{22}/\sqrt{g(\la_i)} \end{array}\right]   \|_F.$$
This leads us to use in our experiments
\begin{equation} \label{eq.estimatorerror}
r(P) := \max_i r(P,\la_i)
\end{equation}
as an estimate for the structured absolute backward error induced by our algorithm, i.e., as an estimate for the numerator of the left-hand side of \eqref{eq.KSRbound}. We emphasize that this is a lower bound for the actual global structured backward error, since it corresponds to a rational matrix that has only one of the computed eigenvalues as a finite zero.

\begin{figure}[!h]
	\begin{center}
		\includegraphics[width=9.5cm, height=7.5cm]{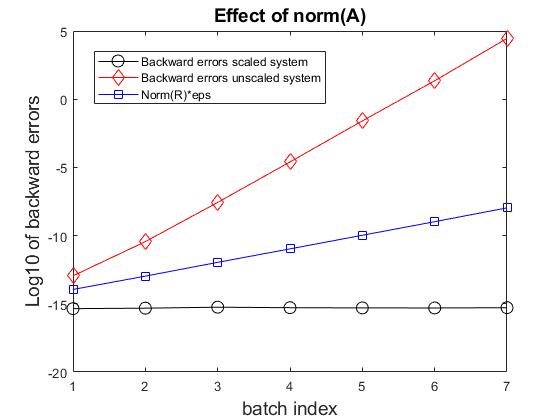}
	\end{center}
	\caption{\label{figure1} Experiment 1: behavior of absolute structured backward errors for increasing values of the norm of $A$.}
\end{figure}

In the first experiment, we investigate the behavior of the structured backward error for rational matrices with matrices $A$ of increasing (large) norms, and with the rest of the matrices in the quadruple $\{A,B,C,D(\la)\}$ having norms of order $1$. The reason why we pay first particular attention to the norm of $A$ is because according to the bound \eqref{eq.KSRbound} the influence of $A$ should be huge because it contributes to $\|S(\la)\|_2$ and also to $f_1$ and $f_3$. For this purpose, we generated with the Matlab function {\tt randn}, 7 batchs of samples of 50 random matrix-tuples $\{A,B,C,D_0,D_1,D_2,D_3\}$, and in each batch indexed with $i$, we multiplied the matrix $A$ by $10^i$, with $i$ going from 1 till 7, in each of the 50 runs of each batch. In each batch, we computed the average of the absolute backward error estimators \eqref{eq.estimatorerror} for both the original matrix-tuples and the scaled ones after applying the procedure in Section \ref{scaling}. In Figure \ref{figure1}, we plot the results of these computations: the horizontal axis represents the index $i$ defining each batch and the vertical axis the logarithm of the average absolute backward errors. Ideally, the absolute backward error should be of order $\epsilon_M \, \|R(\la)\|_F$, where $\epsilon_M$ is the machine precision,  and, so, we also plot this magnitude for the unscaled original data taking in each batch the average of all $\|R(\la)\|_F$ (for the scaled data, this magnitude is always of order $\epsilon_M$ and is not plotted). We observe that the absolute backward errors for the unscaled problem grow very strongly with the index $i$, i.e., with the norm of $A$, and that computing the zeros of a rational matrix by applying the $QZ$ algorithm to the block Kronecker linearization $S(\la)$ is highly structurally backward unstable for large norms of $A$, as predicted by the bound \eqref{eq.KSRbound}. In contrast,
when applying the scaling procedure described in Section \ref{scaling}, this growth is absent and we get perfect structural backward stability for the scaled rational matrix, as predicted by \eqref{eq.KSRbound}.

\begin{figure}[!h]
	\begin{center}
		\includegraphics[width=9.5cm, height=7.5cm]{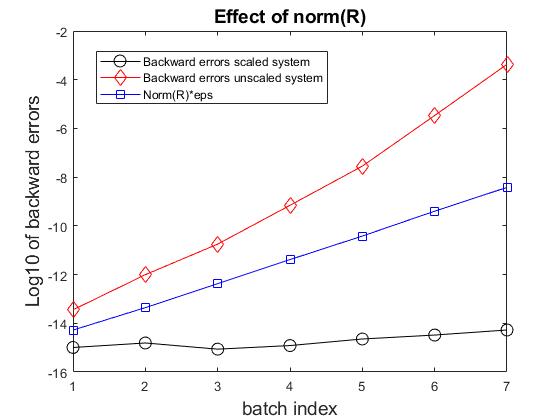}
	\end{center}
	\caption{\label{figure2} Experiment 2: behavior of absolute structured backward errors for increasing values of the norms of $B$, $C$ and $D(\la)$.}
\end{figure}

In the second experiment, we investigate the behavior of the structured backward error for rational matrices with matrices $A$ of norms of order $1$, and with the rest of the matrices in the quadruple $\{A,B,C,D(\la)\}$ having increasing (large) norms. The situation in this experiment is opposite to the one in the first experiment. The matrices are generated following the same pattern of the first experiment except by the fact that once the matrices $\{A,B,C,D_0,D_1,D_2,D_3\}$ are generated with {\tt randn}, $B$ is multiplied by $10^{i/2}$, $C$ by $10^{i/3}$, $D_1$ by $10^{i}$, $D_2$ by $10^{i/2}$ and $D_3$ by $10^{i/3}$, for $i=1, \ldots , 7$. The results are plotted in Figure \ref{figure2} and the conclusions are the same as in the first experiment and are in agreement with our analysis. However, note that the growth of the absolute backward errors of the original unscaled data is much smaller than in the first experiment. This effect is qualitatively expected from the bound \eqref{eq.KSRbound}, since $f_3$ does not depend on the norms of $B$, $C$ and $D(\la)$, but the observed very large quantitative difference is not fully explained by  \eqref{eq.KSRbound}. Possible reasons of this are that, as we have emphasized before, our backward error estimator is a lower bound that may underestimate severely the actual global backward error and/or that the bound in \eqref{eq.KSRbound} overestimates the actual error.

\begin{figure}[!h]
	\begin{center}
		\includegraphics[width=9.5cm, height=7.5cm]{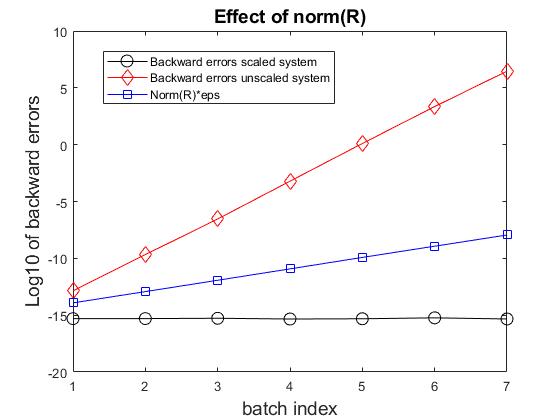}
	\end{center}
	\caption{\label{figure3} Experiment 3:  behavior of absolute structured backward errors for increasing values of the norms of $A$, $B$, $C$ and $D(\la)$.}
\end{figure}

The last experiment we present combines the scalings used in the first and second experiments. That is,  once the matrices $\{A,B,C,D_0,D_1,D_2,D_3\}$ are generated with {\tt randn}, $A$ is multiplied by the factor used in Experiment 1 and $B,C,D_1,D_2$, and $D_3$ are multiplied by the factors used in Experiment 2. Taking into account that the function $f_1$ appearing in the bound \eqref{eq.KSRbound} includes a product of the norm of $A$ times the norm of $B$ and a product of the norm of $A$ times the norm of $C$, we expect backward errors larger than those of Experiment 1. The results are plotted in Figure \ref{figure3}. The errors are indeed larger than those in Figure \ref{figure1}, but just a bit larger. The possible reasons of this small increment of the errors are the same as in the second experiment.

The main conclusion of this section is that our main {\em a priori} structured backward error bound \eqref{eq.KSRbound} identifies correctly the sources of instability of computing the eigenstructure of a rational matrix by applying the $QZ$ algorithm to its block Kronecker linearizations and that the scaling proposed in Section \ref{scaling} leads to structural backward stability.

\section{Conclusions and future work} \label{sec.conclusions}
We have developed the first structured backward error analysis for an algorithm that computes the eigenstructure of a rational matrix. More precisely, the considered algorithm  starts from a rational matrix expressed as in \eqref{eq_rationalmatrix} and computes its eigenstructure by applying a backward stable generalized eigenproblem algorithm to its block Kronecker linearizations described in \eqref{block_kronecker_lin}. As a consequence of this analysis, we have identified the simple sufficient conditions \eqref{eq.sufficient_conditions} for structural backward stability. In the case of rational matrices which do not satisfy these conditions, we have developed an scaling procedure that transforms the original matrix in another one for which structural backward stability is guaranteed. A number of numerical experiments confirming the predictions of the backward error analysis have been performed and discussed. The results in this paper open new research problems in the area of structured backward error analysis, since other representations used in applications of the given rational matrix  should be considered in the future, as well as other families of linearizations.

\appendix
\section{Auxiliary result for Lemma \ref{lem:singval_4}} \label{Appendix}
We prove in this appendix that the matrix
$$ \left[ \begin{array}{c|c} A & B \\ \hline
C & D\end{array} \right]
:= \left[ \begin{array}{c|c} E_k^T \otimes I_{k} & I_{(k+1)} \otimes E_k \\ \hline
F_k^T \otimes I_{k} & I_{(k+1)} \otimes F_k \end{array} \right]
$$
appearing in the proof of Lemma \ref{lem:singval_4} can be transformed by row and column permutations to the direct sum of the following matrices~:
$$  M_1\oplus M_1\oplus M_3 \oplus M_3 \oplus \cdots \oplus M_{2k-1} \oplus  M_{2k-1} \oplus N_{2k},
$$
where the blocks $M_k$ and $N_k$ are as defined in \eqref{MkNk}.
Let us take for example $k=3$, then the matrix looks like
$$ \left[ \begin{array}{c|c} A & B \\ \hline
C & D\end{array} \right]
:= \left[ \begin{array}{ccc|ccccc} I_3 & & & E_3 & & & \\ & I_3 & & & E_3 & &  \\  &  & I_3 & & & E_3 & \\ &  &  &  & & & E_3 \\ \hline
& & & F_3 & & & \\ I_3 & & & & F_3 & &  \\  & I_3 &  & & & F_3 & \\ &  & I_3 &  & & & F_3
\end{array} \right].
$$
There are three submatrices $M_1$, $M_3$ and $M_5$ that take elements $a$, $b$, $c$ and $d$ in the respective blocks $A$, $B$, $C$ and $D$, as indicated below
$$M_1=\left[ \begin{array}{c} b \end{array}\right], \quad M_3=\left[ \begin{array}{ccc} b & a & \\ & c & d \\ & &  b \end{array}\right], \quad M_5=\left[ \begin{array}{ccccc} b & a & & & \\ & c & d & & \\ & & b & a & \\ & & & c & d \\ & & & & b \end{array}\right]$$
and they each start with a leading element in one of the $E_3$ blocks. For instance, $M_1=\left[b_{10,13}\right]$, $M_3$ starts with the leading element $b_{7,9}$ in the third $E_3$ block, and $M_5$ starts with the leading element in the second $E_3$ block~:
$$ M_1=\left[b_{10,13}\right],$$ $$  M_3=\left[ \begin{array}{ccc} b_{7,9} & a_{7,7} & \\ & c_{10,7} & d_{10,14} \\ & &  b_{11,14} \end{array}\right], \quad
M_5=\left[ \begin{array}{ccccc} b_{4,5} & a_{4,4} & & &\\ & c_{7,4} & d_{7,10} & & \\ & &  b_{8,10} & a_{8,8} & \\
& & & c_{11,8} & d_{11,15} \\ & & & & b_{12,15} \end{array}\right].$$
Notice that the $\left[\begin{array}{cc}b &a \end{array}\right]$ and $\left[\begin{array}{cc}c &d \end{array}\right]$ pairs have the same row index and that the $\left[\begin{array}{c}a \\ c \end{array}\right]$ and $\left[\begin{array}{c}d \\ b \end{array}\right]$ pairs have the same column index,
which explains the permutation that has to be constructed to extract the matrix. Also the transitions
$$  b_{7,9} \rightarrow b_{11,14},\quad \mathrm{and} \quad b_{4,5} \rightarrow b_{8,10} \rightarrow b_{12,15}
$$
always go down to the next diagonal element in the next $E_3$ block.
In a similar fashion, one finds another set of submatrices  $M_1$, $M_3$ and $M_5$ that take elements $a$, $b$, $c$ and $d$ in the respective blocks $A$, $B$, $C$ and $D$ in a different order, as indicated below
$$M_1=\left[ \begin{array}{c} d \end{array}\right], \quad M_3=\left[ \begin{array}{ccc} d & c & \\ & a & b \\ & &  d \end{array}\right], \quad M_5=\left[ \begin{array}{ccccc} d & c & & & \\ & a & b & & \\ & & d & c & \\ & & & a & b \\ & & & & d \end{array}\right]$$
and they each start with a trailing element in one of the first three $F_3$ blocks. Finally, the remaining matrix $N_6$ takes elements in the blocks
$A$, $B$, $C$ and $D$ in the following order
$$ N_6=\left[ \begin{array}{ccccccc} b & a & & & & & \\ & c & d & & & & \\ & & b & a & & & \\ & & & c & d & & \\ & & & & b & a \\
& & & & & c & d \end{array}\right]$$
and starts with the leading element in the leading $E_3$ block, and ends with the trailing element in the trailing $F_3$ block.

\bibliographystyle{amsplain}

\end{document}